\documentclass[aop,MSNbibl,seceqn,dvips]{arximspdf}
\usepackage{graphics}

\doi{10.1214/10-AOP553}
\volume{39}
\issue{1}
\pubyear{2011}
\firstpage{139}
\lastpage{175}

\makeatletter
\newcommand{\eqref}[1]{(\ref{#1})}
\newtheorem{theorem}{Theorem}[section]
\newtheorem{lemma}[theorem]{Lemma}
\newtheorem{proposition}[theorem]{Proposition}

\newproclaim{rem}[theorem]{Remark}

\newcommand{\ind}{\mathbf{1}}
\newcommand{\R}{\mathbb{R}}
\newcommand{\Z}{\mathbb{Z}}
\newcommand{\N}{\mathbb{N}}
\renewcommand{\tilde}{\widetilde}
\newcommand{\hatt}{\widehat}
\newcommand{\barr}{\overline}

\newcommand{\bP}{\mathbf P}
\newcommand{\dd}{\mathrm d}

\newcommand{\var}{\operatorname{Var}}

\newcommand{\bbL}{\mathbb L}
\newcommand{\bbP}{\mathbb P}
\newcommand{\bbR}{\mathbb R}

\newcommand{\gb}{\beta}
\newcommand{\gga}{\gamma}
\newcommand{\gd}{\delta}
\newcommand{\gep}{\varepsilon}
\newcommand{\gO}{\Omega}
\newcommand{\gL}{\Lambda}

\newcommand{\tf}{\textsc{f}}

\renewcommand{\ll}{\left\langle}

\newcommand{\rr}{\right\rangle}
\makeatother

\begin{document}

\begin{frontmatter}

\title{Influence of spatial correlation \break for directed polymers}
\runtitle{Influence of spatial correlation for directed polymers}

\begin{aug}
\author[A]{\fnms{Hubert} \snm{Lacoin}\thanksref{T1}\corref{}\ead[label=e1]{lacoin@math.jussieu.fr}}
\thankstext{T1}{Supported by ANR Grant POLINTBIO.}
\runauthor{H. Lacoin}
\affiliation{Universit\`a di Roma Tre}
\address[A]{Dipartemento de Matematica\\
Universit degli Studi Roma Tre\\
Largo San Leonardo Murialdo\\
00146 Rome\\ 
Italy\\
\printead{e1}}
\end{aug}

\received{\smonth{12} \syear{2009}}
\revised{\smonth{4} \syear{2010}}

%
\begin{abstract}
In this paper, we study a model of a Brownian polymer in $\R_+\times
\R^d$,
introduced by Rovira and Tindel [\textit{J. Funct. Anal.} \textbf
{222} (2005) 178--201].
Our investigation focuses mainly on the effect of strong spatial
correlation in the environment in
that model in terms of free energy, fluctuation exponent and volume
exponent. In particular,
we prove that under some assumptions, very strong disorder and
superdiffusivity hold
at all temperatures when $d\ge3$ and provide a novel approach to
Petermann's superdiffusivity result in dimension one
[Superdiffusivity of directed polymers in random environment (2000)
Ph.D. thesis].
We also derive results for a Brownian model of pinning in a nonrandom
potential with power-law decay at infinity.
\end{abstract}

%


\end{frontmatter}

\section{Introduction}

\subsection{Motivation and description of the model}
Much progress has been made lately in the understanding of localization
and delocalization phenomena for random polymer models and especially
for a directed polymer in a random environment (see \cite
{cfCSYrev,cfdenH} for reviews on the subject). The directed polymer
in random
environment was first introduced in a discrete setup, where the polymer
is modeled by the graph of a random walk in $\Z^d$ and the polymer
measure is a modification of the law of a simple random walk on $\Z^d$.
Recently, though, there has been much interest in the corresponding
continuous models, involving Brownian motion rather than simple random
walk (see \cite{cfBTV,cfCYcmp,cfRT,cfViens}), or semicontinuous
models (continuous time and discrete space \cite{cfMRT}, discrete time
and continuous space \cite{cfM,cfP}).

The advantage of these continuous or semicontinuous models is that they
allow the use of techniques from stochastic calculus to derive results
in a simple way. Another advantage is that they are a natural framework
in which to study the influence of spatial correlation in the
environment. In this paper, we investigate the influence of slowly vanishing spatial correlation for
the directed polymer in Brownian environment, first introduced by Rovira and Tindel \cite{cfRT}, which we now describe.

Let $(\omega(t,x))_{(t\in\R_+,x\in\R^d)}$ be a real centered
Gaussian field (under the probability law $\bP$) with covariance
function
%
\begin{equation}\label{defQQ}
\bP[\omega(t,x)\omega(s,y)]=:(t\wedge s) Q(x-y),
\end{equation}
where $Q$ is a continuous nonnegative covariance function going to
zero at infinity (here, and throughout the paper, $\bP[f(\omega)]$ denotes
expectation with respect to $\bP$; analogous notation is used for other
probability laws). Informally, the field can be seen as a summation in
time of independent, infinitesimal translation-invariant fields
$\omega(\dd t,x)$ with covariance function $Q(x-y)\,\dd t$. To avoid
normalization, we assume $Q(0)=1$. We define the random Hamiltonian
formally as
\[
H_{\omega,t}(B)=H_t(B):=\int_0^{t}\omega(\dd s,B_s).
\]
For a more precise definition of $H_t$, we refer to \cite{cfBTV},
Section 2,
where a rigorous meaning is given for the above formula.
Notice that with this definition, $(H_t(B))$ is a centered
Gaussian field indexed by the continuous function $B\in C [0,t]$, with
covariance matrix
\[
\bP\bigl[H_t\bigl(B^{(1)}\bigr)H_t\bigl(B^{(2)}\bigr)\bigr
]:=\int_{0}^{t}Q\bigl(B^{(1)}_s-B^{(2)}_s\bigr)\,\dd s.
\]
For most of the purposes of this article, this could be considered as
the definition of $H_t$.

One defines the (random) polymer measure for inverse temperature $\gb$
as a transformation of the Wiener measure $P$ as follows:
\[
\dd\mu^{\gb,\omega}_t(B):=\frac{1}{Z^{\gb,\omega}_t}\exp(\gb
H_t(B))\,\dd P(B),
\]
where $Z^{\gb,\omega}_t$ is the partition function of the model
\[
Z^{\gb,\omega}_t:=P [\exp(\gb H_t)].
\]

The aim of studying a directed polymer is to understand the behavior of
$(B_s)_{s\in[0,t]}$ under $\mu_t$ when $t$ is large for a typical
realization of the environment~$\omega$.

\subsection{Very strong disorder and free energy}
To study some characteristic properties of the system, it is useful to
consider the renormalized partition function
\[
W^{\gb,\omega}_t=W_t:=P\biggl[\exp\biggl(\int_0^{t} \gb\omega
(\dd s,B_s)-\frac{\gb^2}{2}\,\dd s\biggr)\biggr]=\frac{Z^{\gb
,\omega}_t}{\bP Z^{\gb,\omega}_t}.
\]
It can be checked, without much effort, that $W_t$ is a positive
martingale with respect to
\[
\mathcal F_t:=\sigma\{\omega_s, s\le t\}.
\]
Therefore, it converges to a limit $W_{\infty}$. It follows from a
standard argument (checking that the event is in the tail
$\sigma$-algebra) that
\[
\bP\{ W_{\infty}:=0\}\in\{0,1\}.
\]
Bolthausen first had the idea of studying this martingale for polymers
in a discrete setup \cite{cfB}. He used it to prove that when the
transversal dimension $d$ is larger than $3$ and $\gb$ is small
enough, the behavior of the polymer trajectory $B$ is diffusive under
$\mu_t$. The technique has since been improved by Comets and Yoshida
\cite{cfCY} to prove that whenever $W_{\infty}$ is nondegenerate,
diffusivity holds. The argument of \cite{cfCY} can be adapted to our
Brownian case. When $\bP\{ W_{\infty}:=0\}=0,$ we say that
{\sl weak disorder} holds; the situation where  $\bP\{ W_{\infty}:=0\}= 1$ is referred to as {\sl
strong disorder}.

In the Gaussian setup, a partial annealing argument shows that
increasing $\gb$ increases the influence of disorder. Indeed, for any
$t\ge0,$
\begin{eqnarray*}
 W^{\gb+\gb',\omega}_t &=& P\biggl[\exp\biggl(\int_0^{t}
(\gb+\gb')\omega(\dd s,B_s)-\frac{\gb^2+\gb'^2}{2}\,\dd s\biggl)\biggr]\\
&\stackrel{(\mathcal L)}{=}& P\biggl[\exp\biggl(\int_0^{t} \gb
\omega^{(1)}(\dd s,B_s)+\sqrt{\gb'^2+2\gb\gb'}\omega^{(2)} (\dd s,B_s)\\
&&\hspace{165pt}{}-\frac{(\gb+\gb')^2}{2}\,\dd s\biggr)\biggr]\\
&=:& \hatt W^{\omega^{(1)},\omega^{(2)}}_t,
\end{eqnarray*}
where the equality holds in law and $\omega^{(1)}$, $\omega^{(2)}$
are two
independent Gaussian fields distributed like $\omega$ (we denote by
$\bP^{(1)}\otimes\bP^{(2)}$ the associated probability). This is also
valid for $t=\infty$. Averaging with respect to $\omega^{(2)}$ on the
right-hand side gives
\[
\bP^{(2)} \hatt W^{\omega^{(1)},\omega^{(2)}}_t= W^{\gb,\omega^{(1)}}_t.
\]
Moreover, for a given realization of $\omega^{(1)}$, $\bP^{(2)} \hatt
W^{\omega^{(1)},\omega^{(2)}}_{\infty}=0$ implies that
\[
\hatt W^{\omega^{(1)},\omega^{(2)}}_{\infty}=0, \qquad\bP
^{(2)}\mbox{-a.s.}
\]
and, therefore,
\begin{eqnarray*}
\bP\{ W^{\gb+\gb',\omega}_{\infty}=0\}&\ge& \bP^{(1)} \bigl\{
\bP^{(2)} \hatt W^{\omega^{(1)},\omega^{(2)}}_{\infty}=0\bigr\}\\
&=&\bP^{(1)} \bigl\{W^{\gb,\omega^{(1)}}_\infty=0\bigr\}=\bP\{
W^{\gb,\omega}_{\infty}=0\}.
\end{eqnarray*}
As a consequence, there exists a critical value $\gb_c$ separating the
two phases, that is, there exists $\gb_c\in[0,\infty)$ such that
\begin{eqnarray*}
\gb\in(0,\gb_c) \quad& \Rightarrow& \quad\mbox{weak disorder
holds},\\
\gb>\gb_c \quad& \Rightarrow& \quad\mbox{strong disorder holds}.
\end{eqnarray*}

From the physicist's point of view, it is, however, more natural to
have a definition of strong disorder using free energy. The quantity to
consider is the difference between quenched and annealed free energy.

\begin{proposition}
The a.s. limit
%
\begin{equation}\label{defpt}
p(\gb):=\lim_{t\rightarrow\infty}\frac{1}{t}\log W_t=\lim
_{t\rightarrow\infty}\frac{1}{t}\bP[\log W_t]=:\lim_{t\rightarrow
\infty} p_t(\gb)
\end{equation}
exists and is almost surely constant. The function $\gb\mapsto p(\gb)$
is nonpositive and nonincreasing.
\end{proposition}

We can define
\[
\barr{\gb}_c:=\sup\{\gb>0 \mbox{ such that } p(\gb)=0\}.
\]
It is obvious from the definitions that $\gb_c\le\barr\gb_c$.

For a proof of the existence of the limits above and their equality, we
refer to \cite{cfRT}, Lemma 2.4, Proposition 2.6. The nonpositivity
follows from Jensen's inequality
\[
\bP[\log W_t]\le\log\bP[W_t]=0.
\]

It can be shown (for results in the discrete setup, see
\cite{cfCHptrf,cfCSY}) that an exponential decay of $W_t$
corresponds to a significant localization property of the trajectories.
More precisely, under this condition, it can be shown that two paths
chosen independently with law $\mu^{\gb,\omega}_{t}$ tend to spend a
positive fraction of the time in the same neighborhood. For example,
whenever the left-hand side exists [i.e., everywhere except for
perhaps countably many $\gb$, as $p(\gb)+\gb^2/2$ is a convex
function], we have
\[
\frac{\partial p}{\partial\gb}(\gb)
:=-\lim_{t\to\infty}\frac{1}{t\gb}\bP\biggl[\mu_{\gb
,t}^{\otimes2}\biggl(\int_{0}^tQ\bigl(B^{(1)}_s-B^{(2)}_s\bigr
)\biggr)\biggr];
\]
see \cite{cfCHptrf}, Section 7, where this equality is proved for
directed polymers in $\Z_+\times\Z^d$. It has become customary to refer
to this situation as {\sl very strong disorder}.

It is widely expected that the two notions of strong disorder coincide
outside the critical point and that we have $\gb_c=\barr\gb_c$.
However, it remains an open and challenging conjecture. In \cite{cfCV}
and \cite{cfLac}, it has been shown that for \vspace{1pt}the directed polymer in
$\Z_+\times\Z^d$ with $d=1, 2$ and i.i.d. site disorder, very strong
disorder holds at all temperatures, and it was previously well known
that there is a nontrivial phase transition when $d\ge3$. The same is
expected to hold in continuous space and time if the correlation
function $Q$ decays sufficiently fast at infinity.

\subsection{Superdiffusivity}

Another widely studied issue for directed polymers is the
superdiffusivity phenomenon \cite{cfBQS,cfBTV,cfJ,cfM,cfP,cfPiza,cfS}.
As mentioned earlier, in the weak disorder phase, the
trajectory of the polymer conserves all the essential features of the
nondisordered model (i.e., standard Brownian motion). Therefore, if one
looks at a trajectory up to time $t$, the end position of the chain,
the maximal distance to the origin and the typical distance of a point in
the chain to the origin are all of order $t^{1/2}$. This is one of the
features of a \textit{diffusive} behavior. It is believed that in the
\textit{strong disorder} phase, this property is changed, and that the
quantities mentioned earlier are greater than $t^{1/2}$ (the chain
tends to go farther from the origin to reach a more favorable
environment). Physicists have conjectured that there exists a positive
real number $\xi>1/2$ such that, under $\mu_t$ for large $t$,
\[
\max_{s\in[0,t]}\|B_s\|\approx t^{\xi}.
\]
We refer to $\xi$ as the \textit{volume exponent}. It is believed
that in
the strong disorder phase, $\xi$ does not depend on the temperature and
is equal to the exponent of the associated oriented last-passage
percolation model which corresponds to zero temperature (see
\cite{cfPiza}).
Moreover, the volume exponent should be related to the \textit
{fluctuation exponent}, $\chi>0,$ which describes the fluctuation of
$\log Z_t$ around its average and is defined (in an informal way) by
\[
\var_{\bP}\log Z_ t \approx t^{2\chi} \qquad \mbox{for large $t$}.
\]
The two exponents should satisfy the scaling relation
\[
\chi=2\xi-1.
\]
Moreover, in dimension $1$, an additional \textit{hyperscaling} relation,
$\xi=2\chi$, should hold and it is therefore widely believed that
$\xi=2/3$ and $\chi=1/3$. In larger dimensions, there is no consensus
in the physics literature regarding the exponent values.

Superdiffusivity remains, however, a very challenging issue since, in
most cases, the existence of $\xi$ and $\chi$ has not been rigorously
established.

However, some mathematical results have been obtained in various
contexts related to directed polymers and can be translated informally
as inequalities involving $\xi$ and $\chi$.
\begin{longlist}[$\bullet$]
\item[$\bullet$] For undirected first-passage percolation, Newman and Piza proved
that $\xi\le3/4$ in every dimension \cite{cfNP} and, in
collaboration with Licea \cite{cfLNP}, that $\xi\ge3/5$ in dimension
$2$ (corresponding to $d=1$ for a directed polymer), using geometric
arguments.
\item[$\bullet$] Johansson proved \cite{cfJ} that $\xi=2/3$ and
$\chi=1/3$ for last-passage oriented percolation with i.i.d.
exponential variables on $\N\times\Z$ (this corresponds to the
discrete directed polymer with $\gb=\infty$). The method he employed
relies on exact calculation and it is probably difficult to adapt to
other cases.
\item[$\bullet$] For a discrete-time continuous-space directed
polymer model, Petermann \cite{cfP} proved that for $d=1$, $\xi\ge
3/5$. M\'ejane \cite{cfM} proved for the same model that $\xi\le3/4$
in every dimension. The result of Petermann has recently been adapted
for Brownian polymer in Brownian environment by Bezerra, Tindel and
Viens~\cite{cfBTV}.
\item[$\bullet$] Very recently, Balazs, Quastel and
Sepp\"al\"ainen \cite{cfBQS} computed the scaling exponent for the
Hopf--Cole solution of the KPZ/stochastic Burgers equation, a~problem
that can be interpreted as a $(1+1)$-dimensional directed polymer in a
random environment given by space--time white noise. Their result is
coherent with physical predictions, that is, $\chi=1/3$ and $\xi=2/3,$
and may lead to exact results for other models.
\item[$\bullet$] Sepp\"al\"ainen \cite{cfS} proved, for directed polymers with
log-gamma distributed
weight, that $\xi=2/3$. As in Johansson's case, his result relies on
exact calculations that are specific to the particular distribution of
the environment.
\end{longlist}

In addition, the following relations linking $\xi$ and $\chi$ have been
proven to hold in various contexts:
\begin{eqnarray*}
\chi& \ge& 2\xi-1,\\
\chi& \ge& \frac{1-d\xi}{2},\\
\chi& \le& 1/2
\end{eqnarray*}
(see, e.g., \cite{cfCYcmp} in the case of Brownian polymer in
Poissonian environment), leading, for example, to $\chi\ge1/8$ in
dimension $1+1$.

\subsection{Presentation of the main results}

In this paper, we focus (mainly) on the case where $Q$ has power-law
decay [recall \eqref{defQQ}]. Unless otherwise stated, we will consider
that there exists $\theta>0$ such that
%
\begin{equation}\label{Qpold}
Q(x)\asymp\|x\|^{-\theta} \qquad\mbox{as } \|x\| \to\infty,
\end{equation}
where $\|\cdot\|$ denotes the Euclidean norm in $\R^d$. By
$f(x)\asymp
g(x)$ as $\|x\| \to\infty$, we mean that there exist positive
constants $R$ and $c$ such that
%
\begin{equation}\label{machinchose}
c^{-1} f(x)\le g(x) \le c f(x) \qquad\forall x,\ \|x\|\ge R.
\end{equation}
In the sequel, we also write, for functions of one real variable,
$f(t)\asymp g(t)$ as $t \to\infty$ and $f(t)\asymp g(t)$ as $t \to
0_+$ [with definitions similar to \eqref{machinchose}]. In this setup,
we obtain various results concerning free energy, volume exponent and
fluctuation exponent. These results show that when the spatial
correlation decays sufficiently slowly ($d\ge2$, $\theta\le2$ or
$d=1$, $\theta<1$), the essential properties of the system are changed,
even in a spectacular way for $d\ge3,$ where the weak disorder phase
disappears and we can prove superdiffusivity.

\begin{theorem}\label{ddd333}
We have the following characterization of weak/strong disorder regimes:
\begin{longlist}[(iii)]
\item[(i)] if $d\ge3$ and $\theta>2$, then $\barr\gb_c\ge\gb_c>0$;
\item[(ii)] if $d\ge2$ and $\theta<2$, then $\gb_c=\barr\gb_c=0$;
\item[(iii)] $d=1$, $\gb_c=\barr\gb_c=0$ for any value of $\theta$.
\end{longlist}
\end{theorem}

In the cases where $\barr\gb_c=0$, we obtain sharp bounds on both sides
for the free energy.
\begin{theorem}\label{dddd4444}
For $d\ge2$, $\theta<2$ or $d=1$, $\theta<1,$ we have
\[
p(\gb)\asymp-\gb^{4/(2-\theta)}.
\]
For $d=1$, $Q\in\bbL_1(\R)$ (with no other assumption on the decay),
we have
\[
p(\gb)\asymp-\gb^{4}.
\]
\end{theorem}

\begin{rem}
For $d=1$, $\theta>1$, one can see that Theorem \ref{dddd4444} is
identical to \cite{cfLac}, Theorem 1.5,
suggesting that, in this case, the Brownian model is in the same
\textit{universality class}
as the discrete model. One would have the same conclusion $d=2$,
$\theta>2$, suggesting that
the system does not feel the correlation if $Q$ is in $\bbL_1(\R^d)$.
\end{rem}

\begin{rem} In the cases we have left unanswered,
namely $d=2,\theta\ge2$ and $d=3$, $\theta=2$, the technique
used for the two-dimensional discrete case (see \cite{cfLac})
can be adapted to prove that $\barr\gb_c=0$. Since the method is relatively
complicated and very similar to that which is applied in the discrete case,
we do not develop it here. In these cases, $p(\gb)$ decays faster than
any polynomial around zero.
For $d=2$, $\theta>2$ or $d=3$, $\theta=2$, one would expect to have
\[
p(\gb)\asymp-\exp\biggl(-\frac{c}{\gb^2}\biggr),
\]
while for $d=2$, $\theta=2$, one should have
\[
p(\gb)\asymp-\exp\biggl(-\frac{c}{\gb}\biggr).
\]
However, in both cases, one cannot get a lower bound and an upper bound
that match.
\end{rem}

For $d\ge3,$ Theorem \ref{ddd333} ensures that diffusivity holds at
high temperatures when $\theta>2$. We have proved that, on the other
hand, superdiffusivity holds (in every dimension) for $\theta<2$.

\begin{theorem}\label{supdiffff}
When $d\ge2$ and $\theta<2$ or $d= 1$ and $\theta<1$, we have
\[
\lim_{\gep\to0}\operatorname{lim\,inf}\limits_{t\to\infty} \bP
\mu_t^{\gb,\omega}
\Bigl\{\sup_{0\le s\le t} \|B_s\|\ge\gep
t^{3/(4+\theta)}\Bigr\}=1.
\]
For $d=1$, $Q\in\bbL_1(\bbR),$ we have
\[
\lim_{\gep\to0}\operatorname{lim\,inf}\limits_{t\to\infty} \bP
\mu_t^{\gb,\omega} \Bigl\{\sup_{0\le s\le t} \|B_s\|\ge\gep
t^{3/5}\Bigr\}=1.
\]
\end{theorem}

In the development, we will not go into the details of the proof of the
case $Q\in\bbL_1(\bbR)$ as it is very similar to the proofs of the
other cases and, furthermore, because it is just a minor improvement of
the result of \cite{cfBTV}.

\begin{rem}
The argument we use for our proof uses change-of-measure and coupling
arguments instead
of computation on Gaussian covariance matrices as introduced by
Petermann and later
adapted by Bezerra, Tindel and Viens (see \cite{cfBTV,cfP}).
In our view, this makes the computation much clearer. Besides, our
proof is shorter and goes substantially further.
\end{rem}

\begin{rem}
Another polymer model, namely: Brownian polymer
in a Poissonian environment, has been introduced and studied
by Comets and Yoshida \cite{cfCYcmp}. We would like to
stress that our proofs do not rely on the Gaussian nature
of the environment and that superdiffusivity with exponent
$3/5$, as well as very strong disorder in dimension $1$ and $2$
(in dimension $2$ one
needs to adapt the method used in \cite{cfLac}), can also be proven
for this model by using methods developed in the present paper.
We focus on Brownian polymer mainly because it is the natural model to
study the effect of long-range spatial memory.
\end{rem}

On the other hand, the bound of M\'ejane also holds for this model and
so we present a short proof for it.

\begin{proposition}\label{mmeejj}
For any values of $\gb$, $d$ and arbitrary $Q$, and any $\alpha>3/4$,
we have
\[
\lim_{t\to\infty} \bP\mu_t\Bigl\{\max_{s\in[0,t]}\|B_s\| \ge
t^{\alpha}\Bigr\}=0.
\]
\end{proposition}

\begin{rem}
The two previous results can be interpreted as
\[
\frac{3}{4+(\theta\wedge d)} \le\xi\le3/4.
\]
Taking $\theta$ close to zero ensures that the upper bound for $\xi$ is
optimal as a bound which holds for any correlation function $Q$. To get
a better upper bound (e.g., $\xi\le2/3$ in the one-dimensional case),
one would have to use explicitly the lack of correlation in the
environment.
\end{rem}

Our final result concerns the lower bound on the variance of $\log
Z_t$.

\begin{theorem}\label{abcdef}
For any values of $\gb$, $d$ and $Q$ such that \eqref{Qpold} holds, if
$\alpha$ is such that
\[
\lim_{t\to\infty} \bP\mu_t\Bigl\{\max_{s\in[0,t]}\|B_s\|\ge
t^{\alpha}\Bigr\}=0,
\]
then there exists a constant $c$ (depending on $\gb$, $d$ and $Q$) such
that
\[
\var_{\bP} \log Z_t \ge c t^{1-(\theta\wedge d)\alpha}.
\]
In particular, for every $\gep,$ one can find $c$ (depending on $\gb$,
$d$, $Q$ and $\gep$) such that
\[
\var_{\bP} \log Z_t \ge c t^{(4-3\theta)/4-\gep}.
\]
\end{theorem}

The previous result can be informally written as
\[
\chi\ge\frac{1-(\theta\wedge d)\xi}{2}.
\]
The paper is organized as follows.
\begin{longlist}[$\bullet$]
\item[$\bullet$] In Section \ref{Bpinning}, we study a homogeneous pinning model
in order to derive results that will be of use for the study of
$p(\gb)$.
\item[$\bullet$] In Section \ref{FreeN}, we prove all the results
concerning weak and strong disorder and the free energy of the directed
polymer, that is, Theorems \ref{ddd333} and \ref{dddd4444}.
\item[$\bullet$] In Section \ref{flucvol}, we prove all the results concerning volume
exponent and fluctuation exponent, that is, Theorem \ref{supdiffff},
Proposition \ref{mmeejj} and Theorem \ref{abcdef}.
\end{longlist}

\section{Brownian homogeneous pinning in a power-law tailed potential}
\label{Bpinning}

\subsection{The model and presentation of the results}

In this section, we study a deterministic Brownian pinning model. Most
of the results obtained in this section will be used as tools to prove
lower bounds on the free energy for the directed polymer model, but
they are also of interest in their own right. This pinning model was
recently introduced and studied in a paper by Cranston et
al. \cite{cfCKMV} in the case of a smooth and compactly supported
potential---various results were obtained using the tools of
functional analysis. Our main interest here is in potentials with
power-law decay at infinity.

$V$ is either a bounded continuous nonnegative function of $\R^d$ such
that $V(x)$ tends to zero when $x$ goes to infinity or else
$V(x)=\ind_{\{\|x\|\le1\}}$.

We define the energy of a continuous trajectory up to time $t$,
$(B_{s})_{s\in\bbR}$, to be the Hamiltonian
\[
G_t(B):= \int_{0}^t V(B_s)\,\dd s.
\]
We define $\nu_t^{(h)}$, the Gibbs measure associated with that
Hamiltonian and pinning parameter (or inverse temperature) $h\in\R$ to
be
\[
\dd\nu_t^{(h)}(B):=\frac{\exp(h G_t(B))}{Y^{(h)}_t}\,\dd P( B),
\]
where $Y^{(h)}_t$ denotes the partition function
\[
Y^{(h)}_t:=P [\exp(h G_t(B))].
\]

As for the Brownian polymer, the aim is to investigate the typical
behavior of the chain under $\nu_t$ for large $t$. The essential
question for this model is whether or not the pinning potential $hV$
is sufficient to keep the trajectory of the polymer near the origin,
where $V$ takes larger values. It is intuitively clear that for large
$h$, the potential localizes the polymer near the origin [the distance
remains $O(1)$ as $t$ grows] and that $h\le0$ has no chance to
localize the polymer. Therefore, one must determine whether the polymer
is in the localized phase for all positive $h>0$ or if the phase
transition occurs for some critical value of $h_c>0$. This question was
answered in \cite{cfCKMV} in the case of compactly supported smooth
$V$, where it was shown that localization holds for all $h>0$ if and
only if $d\le2$ (in fact, this issue is strongly related to
recurrence/transience of the Brownian motion).

Answering this question relies on studying the free energy.

\begin{proposition}\label{th:freeee}
The limit
\[
\tf(h):=\lim_{t\to\infty} \frac{1}{t}\log Y^{(h)}_t
\]
exists and is nonnegative. $h\mapsto\tf(h)$ is a nondecreasing,
convex function.
\end{proposition}

We call $\tf(h)$ the \textit{free energy} of the model. We define
\[
h_c=h_c(V):=\inf\{h\dvtx\tf(h)>0\}\ge0.
\]

The existence of the limit is not straightforward. Cranston et al. proved
it in \cite{cfCKMV}, Section 7, in the case of $C^{\infty}$
compactly supported $V$; we will adapt their proof to our case. To
understand why $\tf(h)>0$ corresponds to the localized phase, we note
the following point: convexity allows us to interchange limit and
derivative; therefore, at points where $\tf$ has a derivative,
\[
\frac{\partial\tf}{\partial h}(h)=\lim_{t\to\infty} \frac
{1}{t}\nu_t\biggl[\int_0^t V(B_s)\,\dd s\biggr].
\]

Cranston et al. proved that when $V$ is smooth and compactly supported,
$h_c(V)=0$ for $d=1,2$ and $h_c(V)>0$ when $d\ge3$, with some estimate
on the free energy around $h_c$. We want to see how this result can be
modified when $V$ has power-law decay at infinity. To do so, we will
use \cite{cfCKMV}, Theorem 6.1, for which we present a simplified
version for $V(x)=\ind_{\{\|x\|\le1\}}$ (which is not smooth, but for
which the result still holds by monotonicity of the free energy).

\begin{theorem}[(From \cite{cfCKMV}, Theorem 6.1)]\label{comp}
Let $V\dvtx\R^d\rightarrow\R_+$ be defined as
\[
V(x)=\ind_{\{\|x\|\le1\}}.
\]
Then, for $d=1$, $2,$ we have $h_c=0$ and, as $h\to0+$,
%
\begin{eqnarray}\label{freeen}
\tf(h)&=& 2h^2\bigl(1+o(1)\bigr) \qquad\mbox{for } d=1,\nonumber\\
[-8pt]\\[-8pt]
\tf(h)&=& \exp\bigl(-h^{-1}\bigl(1+o(1)\bigr)\bigr)\qquad \mbox
{for }
d=2.\nonumber
\end{eqnarray}
For $d\ge3$, we have $h_c>0$.
\end{theorem}

\begin{rem}
Although Theorem \ref{comp} is completely analogous to results
obtained for
discrete random walk pinning on a defect line (see \cite{cfBook}),
the methods
used to prove it are very different. Indeed, the technique used in the
discrete setup uses discrete renewal theory and cannot be adapted here.
\end{rem}

Unless otherwise stated, we assume from now on that $V$ has power-law
decay at infinity, that is, that there exists $\theta>0$ such that
%
\begin{equation}\label{poldecay}
V(x)\asymp\|x\|^{-\theta} \qquad\mbox{as } \|x\| \to\infty
\end{equation}
[the notation $\asymp$ was introduced in \eqref{Qpold}]. We prove that
in dimension $d \ge3$, whether or not $h_c$ is equal to zero depends
only on the exponent $\theta$. Furthermore, when the value of $\theta$
varies, so does the critical exponent, which can take any value in
$(1,\infty)$.

\begin{theorem}\label{pinning}
For $d\ge3,$ we have
%
\begin{equation}\label{asymp0}
\theta> 2\quad\Rightarrow\quad h_c(V)>0
\end{equation}
and when $h>0$ is small enough, we have
\[
\lim_{t\to\infty} Y^{(h)}_t<\infty.
\]
Moreover, when $\theta<2,$
%
\begin{equation}\label{asymp1}
\tf(h)\asymp h^{2/(2-\theta)} \qquad\mbox{as } h\to0_+.
\end{equation}
\end{theorem}

For the lower-dimensional cases ($d=1,2$), it follows from
\cite{cfCKMV}, Theorem 6.1, and monotonicity of $\tf$ in $V$ that
$h_c=0$ for all $\theta$. However, some of the features of Theorem
\ref{pinning} still hold.
\begin{theorem}
For $d=2$, $\theta<2$ and $d=1$, $\theta<1$, we have
%
\begin{equation}\label{asymp2}
\tf(h)\asymp h^{2/(2-\theta)} \qquad\mbox{as } h\to0_+.
\end{equation}
\end{theorem}

For the sake of completeness, we also present the result for the case
$d=1$, ${\theta>1}$. The following is the generalization of a result
proved for compactly supported smooth function in \cite{cfCKMV},
Theorem 6.1.
\begin{proposition}\label{L1case}
For $d=1$ and $V\in\bbL_1(\bbR)$ continuous and nonnegative,
%
\begin{equation}\label{L11}
\tf(h)= \frac{\|V\|_{\bbL_1(\R)}^2}{2}h^2\bigl(1+o(1)\bigr)
\qquad\mbox{as } h\to0_+.
\end{equation}
\end{proposition}

In the case $d=2$, $\theta\ge2$, it can be checked that monotonicity
in $V$ of the
free energy and \eqref{freeen} together imply that $h_c>0$ and that
$\tf(h)$ decays
faster that any power of $h$ around $h=0$. In the case $d=1$, $\theta
=1$, monotonicity again implies that $\tf(h)$ behaves like
$h^{2+o(1)}$ around $0+$.

\begin{rem}
In \cite{cfCKMV}, the critical behavior of the free energy is
computed (in a sharp way)
for every dimension, even in the case where $h_c>0$ ($d\ge3$), as it
is done in the discrete case.
As will be seen in the proof, the method used to obtain critical
exponents in the present paper fails
to give any result when $h_c>0$. However, it would be natural to think
that for some value of $\theta>2$,
the critical exponent depends on $\theta$ and is different from the
one obtained in \cite{cfCKMV}, Theorem~6.1.
\end{rem}

\subsection{Preparatory proof}
In this subsection, we give the proofs of results that are easy
consequences of results in \cite{cfCKMV}: existence of the free
energy and an estimate of
the free energy in the case $V\in\bbL^{1}(\bbR)$.

\begin{pf*}{Proof of Proposition \ref{th:freeee}}
Given $V$ and $\gep$, one can find a compactly supported $C^{\infty}$
function $\check V$ such that
\[
\check V \le V \le\check V+ \gep.
\]

We write $\check Y_t^{(h)}$ for the partition function corresponding to
$\check V$. Trivially, we have, for every $t$,
\[
\frac{1}{t}\log\check Y_t^{(h)}\le \frac{1}{t}\log Y_t^{(h)}\le
\frac{1}{t}\log\check Y_t^{(h)}+h\gep.
\]
As proved in \cite{cfCKMV}, Section 7, $\log\check Y_t$ converges as
$t$ goes to infinity so that
\[
\operatorname{lim\,sup}\frac{1}{t}\log Y_t^{(h)}- \operatorname
{lim\,inf}\frac{1}{t}\log Y_t^{(h)}\le h\gep.
\]
The proof can also be adapted to prove the existence of the free energy
for the potential
\[
V(x):=\ind_{\{\|x\|\le1\}}.
\]
We omit the details here.
\end{pf*}

\begin{pf*}{Proof of proposition \ref{L1case}}
First, we prove the upper bound. By the occupation times formula (see,
e.g.,~\cite{cfRY}, page 224), we have
\[
\int_0^t V(B_s)\,\dd s= \int_{\R} L_t^xV(x)\,\dd x,
\]
where $L_t^x$ is the local time of the Brownian motion in $x$ at time
$t$. By Jensen's inequality, we then have
\[
\exp\biggl(h \int_{\R} L_t(x)V(x)\,\dd x\biggr)\le\int_{\R
}\frac{V(x)\,\dd x}{\|V\|_{\bbL_1(\R)}}\exp\bigl(h \|V\|_{\bbL
_1(\R)}L_t^x\bigr).
\]
Moreover, under the Wiener measure with initial condition zero, $L_t^x$
is sto\-chastically bounded (from above) by $L_t^0$ for all $x$ so that
%
\begin{eqnarray}\label{boundL1}
Y_t^{(h)}&\le& \int_{\R}\frac{V(x)\,\dd x}{\|V\|_{\bbL_1(\R)}}
P\bigl[\exp\bigl(h\|V\|_{\bbL_1(\R)}L_t^x\bigr)\bigr]\nonumber
\\ [-8pt]\\ [-8pt]
&\le& P\bigl[\exp\bigl(h \|V\|_{\bbL_1(\R)}L^0_t\bigr)\bigr]\le2
\exp\biggl(\frac{t h^2 \|V\|^2_{\bbL_1(\R)}}{2}\biggr).\nonumber
\end{eqnarray}
In the last inequality, we have used the fact that
$L_t^0\stackrel{(\mathcal{L})}{=}\sqrt{t}|\mathcal N(0,1)|$. Taking the
limit as $t$ tends to infinity gives the upper bound. For the lower
bound, the assumption we have on $V$ guarantees that, given $\gep>0$,
we can find $\check V,$ smooth and compactly supported, such that
\begin{eqnarray*}
\check V &\le& V,\\
\|\check V\|_{\bbL_1(\R)} &\ge& \| V \|_{\bbL_1(\R)}-\gep.
\end{eqnarray*}
Let $\check\tf$ be the free energy associated with $\check V$. By
\cite{cfCKMV}, Theorem 6.1, and monotonicity, we have that for $h$
small enough (how small depends on $\gep$),
\[
\tf(h)\ge\check\tf(h)\ge\biggl(\frac{\|\check V\|_{\bbL_1(\R
)}^2}{2}-\gep\biggr)h^2.
\]
Since $\gep$ was chosen arbitrarily, this gives the lower bound.
\end{pf*}

\subsection{Proof of upper bound results on the free energy}

In this subsection, we prove the upper bounds corresponding to
\eqref{asymp0}, \eqref{asymp1} and \eqref{asymp2}. We will give a brief
sketch of the proof. To do so, we will repeatedly use the following
result.

\begin{lemma}\label{lemrelou}
Let $(a_n)_{n\in\N}$ be a sequence of positive real numbers and
$(p_n)_{n\in\N}$ a sequence of strictly positive real numbers
satisfying $\sum_{n\in\N} p_n=1$.
We then have
%
\begin{equation}\label{prood}
\prod_{n\in\N} a_n\le\sum_{n\in\N} p_na_n^{1/p_n}
\end{equation}
provided the left-hand side is defined. This formula is also valid for
a product with finitely many terms.
\end{lemma}
\begin{pf}
Let $X$ be the random variable whose distribution $P$ is defined as
follows:
\[
P\{X=x\}=\sum_{\{n: (\log a_n)/p_n=x\}}p_n.
\]
The formula considered is just Jensen's inequality:
\[
\exp(P[X])\le P \exp(X).
\]
\upqed
\end{pf}

We now proceed to the proof of delocalization at high temperature for
$\theta>2$ and $d\ge3$. The strategy is to use Lemma \ref{lemrelou} to
bound the partition function by a convex combination of countably many
partition functions of pinning systems with compactly supported
potential and then to use rescaling of Brownian motion to get that all
of these partition functions are uniformly bounded as $t\to\infty$.

\begin{pf*}{Proof of \eqref{asymp0}}
Let $\theta>2$ and $\gep>0$. We define
%
\begin{eqnarray}\label{salusalu}
\barr{V}(x)&:=&\sum_{n=0}^{\infty} \ind_{\{\|x\| \le
2^n\}}2^{-n\theta}\nonumber\\ [-8pt]\\ [-8pt]
&=&\frac{1}{1-2^{-\theta}}\Biggl[\ind_{\{\|x\|\le1\}}+\sum
_{n=1}^{\infty}\ind_{\{2^{n-1}< \|x\|\le
2^{n}\}}2^{-n\theta}\Biggr].\nonumber
\end{eqnarray}
Assuming that $V>0$, it follows from assumption \eqref{poldecay} that
there exist constants $c_1$ and $C_1$ such that
%
\begin{equation}\label{barVV}
c_1 \barr V(x) \le V(x)\le C_1 \barr V (x).
\end{equation}
The proof would also work for $V\ge0$, by defining $\barr V$
differently, that is, dropping the $\ind_{\{\|x\|\le1\}}$ term and
starting the sum from some large $n_0$ at the left-hand side of
\eqref{salusalu}. We restrict our attention to the case $V>0$ here,
only for notational reasons: sums on balls are cleaner (to write) than
sums on annuli.

Hence, for any $p\in(0,1)$ and $h>0$, we have
\begin{eqnarray*}
 Y_t^{(h)}&\le& P\biggl[\exp\biggl(hC_1 \int_{0}^t \barr{V}(B_s)\,
\dd s\biggl)\biggl]\\
&\le& (1-p)\sum_{n=0}^{\infty}p^{n} P\biggl[ \exp\biggl
((1-p)^{-1}p^{-n}hC_1\int_0^t\ind_{\{\|B_s\|
\le2^n\}}2^{-n\theta}\,\dd s \biggr)\biggr]\\
&= &(1-p)\sum_{n=0}^\infty p^n P\biggl[\exp\biggl
((1-p)^{-1}p^{-n}2^{n(2-\theta)}h C_1\int_0^{2^{-2n}t}\ind_{\{\|
B_s\| \le1\}}\,\dd s\biggr)\biggr],
\end{eqnarray*}
where the second inequality uses Lemma \ref{lemrelou} with $p_n:=
(1-p)p^n$ and the last equality is just a rescaling of the Brownian
motion. We choose $p$ such that $p2^{\theta-2}=1$. We get
\[
Y_t^{(h)}\le P\biggl[\exp\biggl((1-p)^{-1}h C_1\int_0^{t}\ind_{\{\|
B_s\|\le1\}}\,\dd s\biggr)\biggr].
\]
For $h$ small enough, Theorem \ref{comp} for $d\ge3$ allows us to
conclude. Moreover (see \cite{cfCYcmp}, Proposition 4.2.1), in this
case, we have
\[
\lim_{t\rightarrow\infty} Y_t^{(h)}:=Y^{(h)}_{\infty}<\infty.
\]
\upqed
\end{pf*}

To prove an upper bound on the free energy when $\theta<2$, we start by
putting aside the contribution given by $V(x)$ when $x$ is large. This
way, we just have to estimate the partition function associated with a
compactly supported potential. In dimension $d\ge3,$ what was done in
the preceding proof will be sufficient to obtain the result. For $d=1$
and $d=2$, we will have to make good use of Theorem \ref{comp}.

\begin{pf*}{Proof of the upper bounds in \eqref{asymp1} and
\eqref{asymp2}}
We start with the case $d=1$, $\theta<1$. From assumption
\eqref{poldecay}, there exists a constant $C_2$ such that for any
$h\le
1$,
%
\begin{eqnarray}\label{boundL12}
V(x)&\le& C_2 h^{\theta/(2-\theta)} \qquad\forall x,\ |x|\ge
h^{-1/(2-\theta)},\nonumber\\ [-8pt]\\ [-8pt]
\int_{|x|\le h^{-1/(2-\theta)}}V(x)\,\dd x &\le&
C_2h^{-(1-\theta)/(2-\theta)}.\nonumber
\end{eqnarray}
We write $\check V(x):= V(x)\ind_{\{|x|\le h^{-1/(2-\theta)}\}}$. We have
\[
V(B_s)\le\check V (B_s)+C_2 h^{\theta/(2-\theta)}
\]
so that
\[
\log Y_t^{(h)} \le t C_2 h^{2/(2-\theta)}+\log P\biggl[\exp\biggl
(h\int_0^t\check V(B_s)\,\dd s\biggr)\biggr].
\]
We know from \eqref{boundL1} and \eqref{boundL12} that the second term
is smaller than
\[
\log2+ t\frac{\|\check V\|^2_{\bbL_1(\R)}}{2}\le\log2+ t C_2^2
h^{2/(2-\theta)},
\]
which is the desired result.

Now, we consider the case $d\ge2$, $\theta<2$. Define $\barr n=\barr
n_{h}:=\lceil|\log h|/[(2-\theta)\log2] \rceil- K$ for some large
integer $K$. We have
\[
\sum_{n>\barr n} 2^{-n\theta} \ind_{\{\|B_s\|\le2^{n}\}}\le \sum
_{n>\barr n} 2^{-n\theta}=\frac{2^{-(\barr n+1)\theta}}{1-2^{-\theta}}.
\]
Therefore [using formula \eqref{barVV}], we can find a constant $C_3$
(depending on $K$ and $C_1$) such that
\begin{eqnarray*}
Y_t^{(h)}&\le& P\biggl[\exp\biggl(C_1 h\int_{0}^{t} \barr V(B_s)\,
\dd s\biggr)\biggr]\\
&\le& e^{C_3t h^{2/(2-\theta)}}P\Biggl[\exp\Biggl(C_1h\int_{0}^{t}
\sum_{n=0}^{\barr n} 2^{-n\theta} \ind_{\{\|B_s\|\le2^{n}\}}\,\dd
s\Biggr)\Biggr],
\end{eqnarray*}
so
%
\begin{eqnarray} \label{stepep}
\frac{1}{t}\log Y_t^{(h)}\le C_3 h^{2/(2-\theta)}+\frac{1}{t}\log P\Biggl
[\exp\Biggl(C_1h\int_{0}^{t} \sum_{n=0}^{\barr n} 2^{-n\theta}
\ind_{\{\|B_s\|\le2^{n}\}}\,\dd s\Biggr)\Biggr].\hspace{-40pt}
\end{eqnarray}

We now have to check that the second term in \eqref{stepep} does not
give a contribution larger than $h^{2/(2-\theta)}$. For any $p$
[we choose $p=2^{(\theta/2)-1}$], as a consequence of Lemma
\ref{lemrelou}, we have
\begin{eqnarray*}
&&P \Biggl[\exp\Biggl(C_1h\int_{0}^{t} \sum_{n=0}^{\barr n}
2^{-n\theta} \ind_{\{\|B_s\|\le2^{n}\}}\,\dd s\Biggl)\Biggr]\\
&&\qquad=P \Biggl[\exp\Biggl(C_1h\int_{0}^{t} \sum_{n=0}^{\barr
n} 2^{-(\barr n-n)\theta} \ind_{\{\|B_s\|\le2^{\barr n-n}\}}\,\dd
s\Biggr)\Biggr]\\
&&\qquad\le\sum_{n=0}^{\barr n}\frac{(1-p)p^n}{1-p^{\barr n+1}} P
\biggl[\exp\biggl(C_1h \int_{0}^{t}\frac{p^{-n}(1-p^{\barr
n+1})}{1-p} 2^{(n-\barr n) \theta} \ind_{\{\|B_s\|\le2^{\barr n-n}\}
}\,\dd s \biggr)\biggr]\\
&&\qquad= \sum_{n=0}^{\barr n}\frac{(1-p)p^n}{1-p^{\barr n+1}}\\
&&\qquad\quad\hphantom{\sum_{n=0}^{\barr n}}{}\times P\biggl[ \exp\biggl(C_1h \int_{0}^{t2^{-2(\barr n-n)}}\frac
{p^{-n}(1-p^{\barr n+1})}{1-p} 2^{(\barr n-n)(2- \theta)} \ind_{\{\|
B_s\|\le1\}}\,\dd s \biggr)\biggr],
\end{eqnarray*}
where the last line is obtained simply by rescaling the Brownian
motion in the expectation. Now, observe that for any $\gep>0$, one can
find a value of $K$ such that
\[
C_1 h\le(1-p)\gep2^{-\barr n(2-\theta)},
\]
so
\[
C_1h\frac{p^{-n}(1-p^{\barr n+1})}{1-p} 2^{(\barr n-n)(2- \theta)}\le
\gep2^{n((\theta/2)-1)}.
\]
Therefore, we have
\begin{eqnarray*}
&&P\Biggl[\exp\Biggl(C_1h\int_{0}^{t} \sum_{n=0}^{\barr n}
2^{-n\theta} \ind_{\{\|B_s\|\le2^{n}\}}\,\dd s\Biggr)\Biggr]\\
&&\qquad\le\max_{n\in\{0,\dots, \barr n\}} P\biggl[\exp\biggl
(\int_{0}^{t2^{-2(\barr n-n)}}\gep2^{n(\theta/2-1)}\ind_{\{\|B_s\|
\le1\}}\,\dd s \biggr)\biggr].
\end{eqnarray*}
For $d\ge3,$ the right-hand side is less than
\[
P\biggl[\exp\biggl(\int_{0}^{t}\gep\ind_{\{\|B_s\|\le1\}}\,\dd
s\biggr)\biggr],
\]
which stays bounded as $t$ goes to infinity. For $d=2$, if $t$ is
sufficiently large and $\gep$ small enough, then Theorem \ref{comp}
allows us to write, for all $n\in\{0,\dots,\barr n\},$
\begin{eqnarray*}
&&\log P\biggl[\exp\biggl(\int_{0}^{t2^{-2(\barr n-n)}}\gep
2^{n((\theta/2)-1)}\ind_{\{\|B_s\|\le1\}}\,\dd s \biggr)\biggr]\\
&&\qquad\le t 2^{-2\barr n} 2^{2n} \exp\bigl(-2\gep^{-1}
2^{n(1-\theta/2)}\bigr).
\end{eqnarray*}
The maximum over $n$ of the right-hand side is attained for $n=0$.
Therefore,
\[
\log P \Biggl[\exp\Biggl(C_1h\int_{0}^{t} \sum_{n=0}^{\barr n}
2^{-n\theta} \ind_{\{\|B_s\|\le2^{n}\}}\,\dd s\Biggr)\Biggr]\le
t2^{-2\barr n}.
\]
Inserting this into \eqref{stepep} ends the proof.
\end{pf*}

\subsection{Proof of lower bounds on the free energy}

In this section, we prove the lower bounds for the asymptotics
\eqref{asymp1} and \eqref{asymp2}. This is substantially easier than
what has been done for upper bounds. Here, one just needs to find a
compactly supported potential which is bounded from above by $V$ and
that gives the appropriate contribution.

For any $n\in\N$,
\[
Y_t^{(h)}\ge P \biggl[ \exp\biggl( c_1 h \int_0^t \barr V(B_s)\,\dd
s\biggr)\biggr]
\ge P\biggl[\exp\biggl(c_1 h 2^{-n\theta}\int_{0}^t \ind_{\{ \|
B_s\|\le2^{-n}\}}\,\dd s\biggr)\biggr].
\]
Rescaling the Brownian motion, we get
\[
Y_t^{(h)}\ge P\biggl[\exp\biggl(c_1 h 2^{n(2-\theta)}\int_{0}^{t2^{-2n}}
\ind_{\{ \|B_s\|\le1\}}\,\dd s\biggr)\biggr].
\]
We can choose $n=n_{h}=\lceil|\log h|/[(2-\theta)\log2] \rceil+ K$
for some integer $K$. Let $C_4>0$ be such that
\[
\lim_{t\rightarrow\infty} \frac{1}{t}\log P \biggl[\exp\biggl
(\int_{0}^tC_4 \ind_{\{\|B_s\|\le1\}}\,\dd s\biggr)\biggr] \ge1.
\]
By choosing $K$ large enough, we can get
\[
\log Y_t^{(h)}\ge \log P \biggl[\exp\biggl(C_4 \int
_{0}^{t2^{-2n_{h}}} \ind_{\{ \|B_s\|\le1\}}\,\dd s\biggr)\biggr]
\ge t 2^{-2n_{h}}.
\]
From this, we get that $h_c=0$ and that
\[
\hspace*{118pt}\tf(h)\ge2^{-2(K+1)} h^{2/(2-\theta)}.\hspace*{118pt}\qed
\]

\begin{rem}
The above proofs indicate that under the measure $\nu_t{(h),}$ when $t$
is large and $\theta<2$, $d\ge2$ or $\theta<1$, $d=1$, the typical
distance of the polymer chain $(B_s)_{s\in[0,t]}$ to the origin is of
order $h^{-1/(2-\theta)}$.
\end{rem}

\section{The directed polymer free energy} \label{FreeN}

\subsection{Lower bounds on the free energy, the second moment method
and replica coupling}
In this section, we make use of the result obtained for homogeneous
pinning models to get some lower bounds on the free energy and prove
the corresponding halves of Theorems \ref{ddd333} and \ref{dddd4444}.
The partition function of a homogeneous pinning model appears naturally
when one computes the second moment of $W_t$.

We start with a short proof of the fact that weak disorder holds for
small $\gb$ if $d\ge3$, $\theta>2$.

\begin{pf}
It is sufficient to show that $W_t$ converges in $L_2$ for $\gb$
sufficiently small. We have
\begin{eqnarray*}
\bP[ W_t^2]&=& \bP\biggl[ P^{\otimes2} \exp\biggl(\int
_0^{t}\bigl[ \gb\omega\bigl(\dd s,B^{(1)}_s\bigr)+ \gb\omega
\bigl(\dd s,B^{(2)}_s\bigr)\bigr]-\gb^2\,\dd s\biggl)\biggr]\\
&=& P^{\otimes2} \biggl[\exp\biggl(\gb^2\int_{0}^t Q\bigl
(B_s^{(1)}-B_s^{(2)}\bigr)\,\dd s\biggr)\biggr].
\end{eqnarray*}
The left-hand side is the partition function of the homogeneous pinning
model described in the first section. Therefore, the result is a simple
consequence of Theorem \ref{pinning}.
\end{pf}

We now prove the lower bound on the free energy corresponding to
Theorem~\ref{dddd4444}. We use a method called \textit{replica
coupling}. The idea of using such a method for directed polymers came
in \cite{cfLac} and was inspired by a work on the pinning model of
Toninelli \cite{cfTcmp}.

\begin{pf}
Define, for $\gb>0$, $r\in[0,1],$
\[
\Phi_t(r,\gb):=\frac{1}{t}\bP\biggl[\log P\exp\biggl(\int
_0^{t}\sqrt{r}\gb\omega(\dd s,B_s)-r\gb^2/2\,\dd s\biggr)\biggr]
\]
and for $\gb>0$, $r\in[0,1]$, $\lambda>0,$
\begin{eqnarray*}
\Psi_t(r,\lambda,\gb)&:=&\frac{1}{2t}\bP\biggl[ \log P^{\otimes
2} \exp\biggl(\int_0^{t}\sqrt{r}\gb\bigl[\omega\bigl(\dd
s,B^{(1)}_s\bigr)+ \omega\bigl(\dd s,B^{(2)}_s\bigr)\bigr]\\
&&\hspace*{105pt}{}+\gb^2\bigl[\lambda Q\bigl
(B_s^{(1)}-B_s^{(2)}\bigr)-r\bigr]\,\dd s\biggr)\biggr]\\
&=:&\frac{1}{2t}\bP\bigl[ \log
P^{\otimes2} \exp\bigl(\hatt H_t\bigl(B^{(1)},B^{(2)},r,\lambda
\bigr)\bigr)\bigr].
\end{eqnarray*}
The function $r\mapsto\Phi_t(r,\gb)$ satisfies [recall the
definition of $p_t$ in \eqref{defpt}]
\[
\Phi_t(0,\gb)=0 \quad\mbox{and} \quad\Phi_t(1,\gb)=p_t(\gb).
\]
In the sequel, we use the following version of the Gaussian integration
by parts formula. The proof is straightforward.
\begin{lemma}
Let $(\omega_1,\omega_2)$ be a centered two-dimensional Gaussian
vector. We
have
\[
\bP[\omega_1f(\omega_2)]:=\bP[\omega_1\omega_2]\bP[f'(\omega_2)].
\]
\end{lemma}

Using this formula we get that
%
\begin{equation}\label{phhhi}
\frac{\dd}{\dd r}\Phi_t(r,\gb)=-\frac{\gb^2}{2t}\bP\biggl[\bigl
(\mu^{(\sqrt{r}\gb)}_t\bigr)^{\otimes2}
\biggl[\int_0^t Q\bigr(B_s^{(1)}-B_s^{(2)}\bigr)\,\dd s\biggr
]\biggr].
\end{equation}
Doing the same for $\Psi_t$, we get
%
\begin{eqnarray}\label{psssi}
&&\frac{\dd}{\dd r} \Psi_t(r,\lambda,\gb)\nonumber\\
&&\qquad=\frac{\gb^2}{2t}\bP\biggl[\frac{P^{\otimes2}e^{\hatt
H_t(B^{(1)},B^{(2)},r,\lambda)} \int_0^t Q(B_s^{(1)}-B_s^{(2)})\,\dd
s}{P^{\otimes2}e^{\hatt
H_t(B^{(1)},B^{(2)},r,\lambda)}}\biggr]\nonumber\\ [-8pt]\\ [-8pt]
&&\qquad\quad{}-\frac{\gb^2}{t}\bP\biggl[\frac{P^{\otimes
4}e^{\hatt H_t(B^{(1)},B^{(2)},r,\lambda)+\hatt H_t(B^{(3)},
B^{(4)},r,\lambda)}\int_{0}^{t}Q(B_s^{(1)}-B_s^{(3)})\,\dd
s}{P^{\otimes4}e^{\hatt H_t(B^{(1)},B^{(2)},r,\lambda)+\hatt
H_t(B^{(3)}, B^{(4)},r,\lambda)}}\biggr]\nonumber\\
&&\qquad\le\frac{\gb^2}{2t}\bP\biggl[\frac{P^{\otimes2}e^{\hatt
H_t(B^{(1)},B^{(2)},r,\lambda)}\int_0^t Q(B_s^{(1)}-B_s^{(2)})\,\dd
s}{P^{\otimes2}e^{\hatt H_t(B^{(1)},B^{(2)},r,\lambda)}}\biggr
]=\frac{\dd}{\dd\lambda}
\Psi_t(r,\lambda,\gb).\nonumber
\end{eqnarray}
This implies that for every $r\in[0,1],$ we have
\[
\Psi_t(r,\lambda,\gb)\le\Psi(0,\lambda+r,\gb).
\]
In view of \eqref{phhhi} and \eqref{psssi}, using convexity and
monotonicity of $\Psi_t$ with respect to $\lambda$ and
$\Psi_t(r,0,\gb)=\Phi_t(r,\gb),$ we have
\begin{eqnarray*}
-\frac{\dd}{\dd r}\Phi_t(r,\gb)&=&\frac{\dd}{\dd\lambda}\Psi
_t(r,\lambda,\gb)\bigg|_{\lambda=0}\\
&\le&\frac{\Psi_t(r,2-r,\gb)-\Phi_t(r,\gb)}{2-r}\le
\Psi_t(0,2,\gb)-\Phi_t(r,\gb),
\end{eqnarray*}
where the last inequality uses the fact that $r\le1$. Integrating this
inequality between zero and one, we get
%
\begin{equation}\label{troidix}
p_t(\gb)\ge(1-e)\Psi_t(0,2,\gb),
\end{equation}
where
\begin{eqnarray*}
\Psi_t(0,2,\gb)&=&\frac{1}{2t}\log P^{\otimes2} \exp\biggl[ 2\gb
^2 \int_0^t Q\bigl(B_s^{(1)}-B_s^{(2)}\bigr)\,\dd s\biggr]\\
&=& \frac{1}{2t}\log P\biggl[ \exp\biggl( 2\gb^2 \int_0^t Q\bigl(\sqrt
{2}B_s\bigr)\,\dd s\biggr)\biggr]=:\frac{1}{2t}\log Y_t.
\end{eqnarray*}
Here, $Y_t$ is the partition function of a homogeneous pinning model
with potential $Q(\sqrt{2}\cdot)$ and pinning parameter $2\gb^2$.
Therefore, we know from Theorem \ref{pinning} and Proposition
\ref{L1case} that
\[
\lim_{t\rightarrow\infty}\frac{1}{t}\log Y_t\asymp\gb
^{4/(2-\theta)}\quad\mbox{or}\quad\lim_{t\rightarrow\infty}\frac
{1}{t}\log Y_t\asymp\gb^4
\]
(where the case to be considered depends on the assumption we have on
$Q$). This, combined with \eqref{troidix}, gives the desired bound.
This completes the proof.
\end{pf}

\subsection[Proof of upper bounds on the free energy (Theorem 1.3)]{Proof of
upper bounds on the free energy (Theorem \protect\ref{dddd4444})}

The technique of the proof is mainly based on a change-of-measure
argument. This method was developed and first usedt for pinning models
\cite{cfGLT} and adapted for directed polymer in \cite{cfLac}. Here,
we have to adapt it to the continuous case and take advantage of the
occurrence of spatial memory in the environment. We briefly sketch an
outline of the proof.
\begin{longlist}[$\bullet$]
\item[$\bullet$] First, we use Jensen's inequality to reduce the proof to
estimating a fractional moment (a noninteger moment) of $W_t$.
\item[$\bullet$] We decompose $W_t$ into different contributions corresponding to
paths along a corridor of fixed width.
\item[$\bullet$] For each corridor, we slightly change the measure via a tilting
procedure which lowers the value of $\omega$ in the corridor.
\item[$\bullet$] We use the change of measure to estimate the fractional moment of
each contribution.
\end{longlist}

We start by stating a trivial lemma which will be useful for our proof
and for the next section.

\begin{lemma}\label{gaussest}
Let $(\omega_x)_{x\in X}$ be a Gaussian field indexed by $X$ defined
on the probability space
$(\Omega,\bbP,\mathcal F)$, closed under linear combination. Define the
measure $\tilde\bbP$ as
\[
\frac{\dd\tilde\bbP}{\dd\bbP}=\exp(\omega_{x_0}-\var\omega_{x_0}/2).
\]
Then, under $\bbP$, $(\omega_x)_{x\in X}$ are still Gaussian variables,
their covariance remain unchanged and their expectation is equal to
\[
\tilde\bbP[\omega_x]= \bbP[\omega_x\omega_{x_0}].
\]
\end{lemma}

We now proceed to the proof. Set $\gamma\in(0,1)$ (in the sequel we
will choose $\gamma=1/2$). We note that
\[
\bP[\log W_t]= \frac{1}{\gga}\bP[\log W_t^{\gga}]\le\frac
{1}{\gga}\log\bP [W_t^{\gga}].
\]
For this reason, we have
%
\begin{equation}\label{momentfrac}
p(\gb)\le\frac{1}{\gga}\liminf_{t\to\infty} \frac{1}{t}\log\bP
W_t^{\gga}.
\end{equation}
Therefore, our aim is to prove that $\bP[ W_t^{\gga}]$ decays
exponentially. Fix $T:=C_1 \gb^{-4/(2-\theta)}$. Choose $t:=Tn$
($n$ is meant to tend to infinity). For $y=(y^1,\dots,\break y^d)\in\Z^d$,
define $I_y:=\prod_{i=1}^d [y^i\sqrt{T},(y^i+1)\sqrt{T})$ (where
$\prod$ here denotes interval product).

We decompose the partition function $W_t$ into different contributions
corresponding to different families of paths. We have
\[
W_t:=\sum_{y_1,\dots,y_n\in\Z^d} W_{(y_1,\dots,y_n)},
\]
where
\[
W_{(y_1,\dots,y_n)}:= P \biggl[\exp\biggl(\int_0^{t} \gb\omega
(\dd s,B_s)-\gb^2/2\,\dd s\biggr) \ind_{\{ B_{kT}\in I_{y_k},
\forall k=1,\ldots,n\}}\biggr].
\]

We use the inequality $(\sum a_i)^{\gga}\le\sum a_i^{\gga}$, which
holds for an arbitrary collection of positive numbers, to get
%
\begin{equation}\label{dec}
\bP[W_t^{\gga}]\le\sum_{y_1,\dots,y_n\in\Z^d} \bP\bigl
[W^{\gamma}_{(y_1,\dots,y_n)}\bigr].
\end{equation}

In order to bound the right-hand side of \eqref{dec}, we use the
following change-of-measure argument: given $Y=(y_1,\dots,y_n)$, and
$\tilde\bP_Y$ a probability measure on $\omega,$ we have
%
\begin{eqnarray}\label{hollder}
\bP W^{\gga}_{(y_1,\dots,y_n)}&=&\tilde\bP_Y \biggl[\frac{\dd\bP
}{\dd\tilde\bP_Y}
W^{\gga}_{(y_1,\dots,y_n)}\biggr]\nonumber\\ [-8pt]\\ [-8pt]
&\le& \biggl(\bP\biggl[\biggl(\frac{\dd\bP}{\dd\tilde\bP
_Y}\biggr)^{\gga/(1-\gga)}\biggr]\biggr)^{1-\gga}\bigl(\tilde\bP
_Y\bigl[ W_{(y_1,\dots,y_n)}\bigr]\bigr)^{\gga}.\nonumber
\end{eqnarray}
One then needs to find a good change of measure to apply this
inequality. Let $C_2$ be a (large) fixed constant. Define the blocks
$A_k$ by
\begin{eqnarray*}
A_k&:=&[(k-1)T,kT]\times\prod_{i=1}^d\bigl[y^{i}_{k-1}-C_2\sqrt
{T},y^{i}_{k-1}+C_2\sqrt{T}\bigr],\\
\barr A_k &:=&\prod_{i=1}^d\bigl[y^{i}_{k-1}-C_2\sqrt
{T},y^{i}_{k-1}+C_2\sqrt{T}\bigr],\\
J_Y&:=&\bigcup_{k=1}^n A_k,
\end{eqnarray*}
with the convention that $y_0=0$. Moreover, we define the random
variable
\begin{eqnarray*}
\gO_k&:=&\frac{\int_{A_k}\omega(\dd t,x)\,\dd x}{\sqrt{T\int
_{\barr A^2_k}Q(x-y)\,\dd x\,\dd y}},\\
\Omega_Y&:=&\sum_{k=1}^n \gO_k.
\end{eqnarray*}
Note that, with this definition, $(\Omega_k)_{k\in\{1,\dots,n\}}$ are
standard centered independent Gaussian variables. Define $\tilde\bP_Y$
by
\[
\frac{\dd\tilde\bP_Y}{\dd\bP}:=\exp(-\Omega_Y-n/2).
\]
From this definition and using the fact that $\gga=1/2$, we have
%
\begin{equation}\label{dete}
\biggl(\bP\biggl[\biggl(\frac{\dd\bP}{\dd\tilde\bP_Y}\biggr
)^{\gga/(1-\gga)}\biggr]\biggr)^{1-\gga}=\exp(n/2).
\end{equation}
We also define the measure $\tilde\bP_1$ by
\[
\frac{\dd\tilde\bP_1}{\dd\bP}:=\exp(-\gO_1-1/2).
\]
We now consider the expectation of $W_{(y_1,\dots,y_n)}$ with respect
to $\tilde\bP_Y$. As the covariance structure of the Gaussian field
remains the same after the change of measure (cf. Lemma
\ref{gaussest}), we have
\begin{eqnarray*}
&&\tilde\bP_Y\bigl[ W_{(y_1,\dots,y_n)}\bigr]= P\exp\biggl(\gb
\tilde\bP_Y\biggl[\int_0^t \omega(\dd s,B_s)\biggr]\biggr)\ind
_{\{ B_{kT}\in I_{y_k}, \forall k=1,\ldots, n\}}\\
&&\hphantom{\tilde\bP_Y\bigl[ W_{(y_1,\dots,y_n)}\bigr]} =
P^{(0)}_{O} \biggl[\exp\biggl( \gb\tilde\bP_Y\biggl[\int_0^T
\omega\bigl( \dd s,B^{(0)}_s\bigr)\biggr]\biggr) \ind_{\{
B^{(0)}_{T}\in I_{y_1}\}}\\
&&\quad\hphantom{\tilde\bP_Y\bigl[ W_{(y_1,\dots,y_n)}\bigr]}
{}\times P^{(1)}_{B_T}\biggl[\exp\biggl(\gb\tilde\bP_Y\biggl[
\int_0^T \omega\bigl( \dd(s+T),B^{(1)}_s\bigr)\biggr]\biggr)\\
&&\quad\hphantom{\tilde\bP_Y\bigl[ W_{(y_1,\dots,y_n)}\bigr]}
{}\times\ind_{\{ B^{(1)}_{T}\in I_{y_2}\}}\cdots\ind_{\{
B^{(n-2)}_{T}\in I_{y_{n-1}}\}}\\
&&\quad\hphantom{\tilde\bP_Y\bigl[ W_{(y_1,\dots,y_n)}\bigr
]}{}\times P^{(n-1)}_{B_{T}}\biggl[\exp\biggl(\gb\tilde\bP
_Y\biggl[\int_0^T \omega\bigl(\dd\bigl(s+(n-1)T\bigr
),B^{(n-1)}_s\bigr)\biggr]\biggr)\\
&&\hspace{184pt}\quad\hphantom{\tilde\bP_Y\bigl[ W_{(y_1,\dots
,y_n)}\bigr]}{}\times\ind_{\{ B^{(n-1)}_{T}\in I_{y_n}\}}\biggr
]\cdots\biggr]\\
&&\hphantom{\tilde\bP_Y\bigl[ W_{(y_1,\dots,y_n)}\bigr]}\le\prod
_{k=1}^{n} \max_{x\in I_{y_{k-1}}} P_x \biggl[\exp\biggl(\gb
\tilde\bP_Y\biggl[\int_0^T \omega\bigl(\dd(s+kT),B_s\bigr
)\biggr]\biggr)\ind_{\{ B_{T}\in I_{y_k}\}}\biggr]\\
&&\hphantom{\tilde\bP_Y\bigl[ W_{(y_1,\dots,y_n)}\bigr]} = \prod
_{k=1}^{n}\max_{x\in I_O} P_x \biggl[\exp\biggl(\gb \tilde\bP
_1\biggl[\int_0^T \omega(\dd s,B_s)\biggr]\biggr)\ind_{\{
B_{T}\in I_{y_k-y_{k-1}}\}}\biggr].
\end{eqnarray*}
The first equality is obtained by the use of the Markov \vspace{1pt}property for
Brownian motion: $P^{(i)}_{B^{(i-1)}_T}$ denotes the Wiener measure for
the Brownian \vspace{1pt} motion conditioned to start at $B^{(i-1)}_T$ and which,
conditionally on $B^{(i-1)}_T$, is independent of all the $P^{(k)}$,
$k<i$; $\omega(\dd(s+kT,B_s)$ denotes the time increment of the
field at
$(s+KT,B_s)$. The inequality is obtained by maximizing over $x\in
I_{y_k}$ for intermediate points; $P_x$ is just the Wiener measure with
initial condition $x$. The last equality just uses translation
invariance. Returning to \eqref{dec}, and using \eqref{hollder} and
\eqref{dete}, we get
%
\begin{equation}\label{machin}
\qquad\bP W_t^{\gga}\le e^{n/2}\biggl[\sum_{y\in\Z^d} \biggl
(\max_{x\in I_O} P_x\biggl[\exp\biggl(\gb \tilde\bP_1\biggl
[\int_0^T \omega(\dd s,B_s)\biggr]\biggr)\ind_{\{ B_{T}\in I_y\}
}\biggr)^{\gga}\biggl]\biggl]^n.
\end{equation}
We are able to prove that the right-hand side decays exponentially with
$n$ if we are able to show that
%
\begin{equation}\label{troidiz}
\sum_{y\in\Z^d} \biggl(\max_{x\in I_O} P_x\biggl[\exp\biggl
(\gb \tilde\bP_1\biggl[\int_0^T \dd\omega(s,B_s)\biggr]\biggr
)\ind_{\{ B_{T}\in I_y\}}\biggr]\biggr)^{\gga}
\end{equation}
is small. To do so, we have to estimate the expectation of the
Hamiltonian under $\tilde\bP_1$. We use Lemma \ref{gaussest} and get
%
\begin{eqnarray}\label{modifexp}
-\tilde\bP_1\biggl[\int_0^T \omega(\dd s,B_s)\biggr]=\bP
\biggl[\gO_1\int_0^T \omega( \dd
s,B_s)\biggr]
=\frac{\int_{A_1} Q(x-B_s)\,\dd x\,\dd s}{\sqrt{T\int_{\barr
A_1^2} Q(x-y)\,\dd x\,\dd
y}}.\hspace*{-50pt}
\end{eqnarray}
The above quantity is always positive. However, it depends on the
trajectory $B$. One can check that, when $d\ge2$, $\theta<2$ or $d=1$,
$\theta<1,$ the assumption of polynomial decay for $Q$ implies that [we
leave the case $d=1$, $Q\in\bbL_1(\R)$ to the reader]
\[
\int_{\barr A_1^2} Q(x-y)\,\dd x\,\dd y\asymp T^{d-\theta/2}.
\]
To control the numerator on the right-hand side of \eqref{modifexp},
we need an assumption on the trajectory. We control the value for
trajectories that stay within $A_1$. For all trajectories
$(s,B_s)_{s\in[0,T]}$ that stay in $A_1,$ we trivially have
\begin{eqnarray*}
&&\int_0^T\int_{[-C_2\sqrt{T},C_2\sqrt{T}]^d}Q(x-B_s)\,\dd x\,\dd
s\\
&&\qquad\ge T\min_{y\in[-C_2\sqrt{T}, C_2\sqrt{T}]^d}\int
_{[-C_2\sqrt{T}, C_2\sqrt{T}]^d}Q(x-y)\,\dd x
\end{eqnarray*}
and the right-hand side is asymptotically equivalent to
$T^{1+(d-\theta)/2}$.

Altogether, we get that there exists a constant $C_3$ (depending on
$C_2$) such that, uniformly on trajectories staying in $A_1$,
%
\begin{equation}\label{esta1}
\tilde\bP_1\biggl[\int_0^T \omega(\dd s,B_s)\biggr]\le
-C_3T^{(2-\theta)/4}.
\end{equation}
The distribution of the Brownian motion allows us to find, for any
$\gep>0$, $R=R_{\gep}$ such that
\begin{eqnarray*}
&&\sum_{\|y\|\ge R} \biggl(\max_{x\in I_O} P_x\biggl[\exp\biggl
(\gb \tilde\bP_1\biggl[\int_0^T \omega(\dd s,B_s)\biggr]\biggr
)\ind_{\{ B_{T}\in I_y\}}\biggr]\biggr)^{\gga}\\
&&\qquad\le\sum_{\|y\|\ge R}\Bigl(\max_{x\in I_O} P_x\{B_T\in I_y\}
\Bigr)^{\gga} \le\gep,
\end{eqnarray*}
where the first inequality simply uses the fact that
$\tilde\bP_1(\cdots)$ is negative. For the terms corresponding to
$\|y\|< R,$ we use the rough bound
\begin{eqnarray*}
&&\sum_{\|y\|<R} \biggl(\max_{x\in I_O} P_x\biggl[\exp\biggl(\gb \tilde
\bP_1\biggl[\int_0^T \omega(\dd s,B_s)\biggr]\biggr)\ind_{\{
B_{T}\in I_y\}}\biggr]\biggr)^{\gga}\\
&&\qquad\le(2R)^d \biggl(\max_{x\in I_O} P_x\biggl[\exp\biggl(\gb
\tilde\bP_1\biggl[\int_0^T \omega(\dd s,B_s)\biggr]\biggr)\biggr
]^{\gga}\biggr).
\end{eqnarray*}
Set $\delta:=(\gep/(2R)^d)^{1/\gga}$. The remaining task in
order to find a good bound on \eqref{troidiz} is to show that
\[
\max_{x\in I_O} P_x\exp\biggl(\gb \tilde\bP_1\biggl[\int_0^T
\omega(\dd s,B_s)\biggr]\biggr)\le\delta.
\]
To get the above inequality, we separate the right-hand side into two
contributions: trajectories that stay within $A_1$ and trajectories
that go out of $A_1$. Bounding these contributions gives
\begin{eqnarray*}
&&\max_{x\in I_O} P_x\biggl[\exp\biggl(\gb \tilde\bP_1\biggl
[\int_0^T \omega(\dd s,B_s)\biggr]\biggr)\biggr]\\
&&\qquad\le P\Bigl\{\max_{s\in[0,T]} |B_s|\ge|C_2-1|\sqrt
{T}\Bigr\}\\
&&\qquad\quad{}+ \max_{x\in I_O} P_x\biggl[\exp\biggl(\gb \tilde
\bP_1\biggl[\int_0^T\omega(\dd s,B_s)\biggr]\biggr)\ind_{\{
(s,B_s):s\in[0,T]\subset A_1\}}\biggr],
\end{eqnarray*}
where the first term in the right-hand side is an upper bound for
\[
\max_{x\in I_O} P_x\biggl[\exp\biggl(\gb \tilde\bP_1\biggl[\int
_0^T \omega(\dd s,B_s)\biggr]\biggr)\ind_{\{(s,B_s):s\in
[0,T]\nsubseteq A_1\}}\biggr].
\]

We can fix $C_2$ so that the first term is less that $\gd/2$. Equation
\eqref{esta1} guaranties that the second term is less than
\[
\exp\bigl(-\gb C_3 T^{(2-\theta)/4}\bigr)=\exp\bigl
(-C_3C_1^{(2-\theta)/4}\bigr)\le\gd/2,
\]
where the last inequality is obtained by choosing $C_1$ sufficiently
large. We have thus shown that \eqref{troidiz} is less than $2\gep$.
Combining this result with \eqref{machin} and \eqref{momentfrac}
implies (for $\gep$ small enough) that
\[
p(\gb)\le-\frac{1}{T},
\]
which is the desired result. \qed

\section{Fluctuation exponent and volume exponent}\label{flucvol}

In this section, we prove Theorem \ref{supdiffff}, Proposition
\ref{mmeejj} and Theorem \ref{abcdef}. We give a sketch of our proof
for the superdiffusivity result in dimension one. The idea is to
compare the energy gain and the entropy cost for going to a distance $t^{\alpha}$
away from the origin.

\begin{figure}[b]

\includegraphics{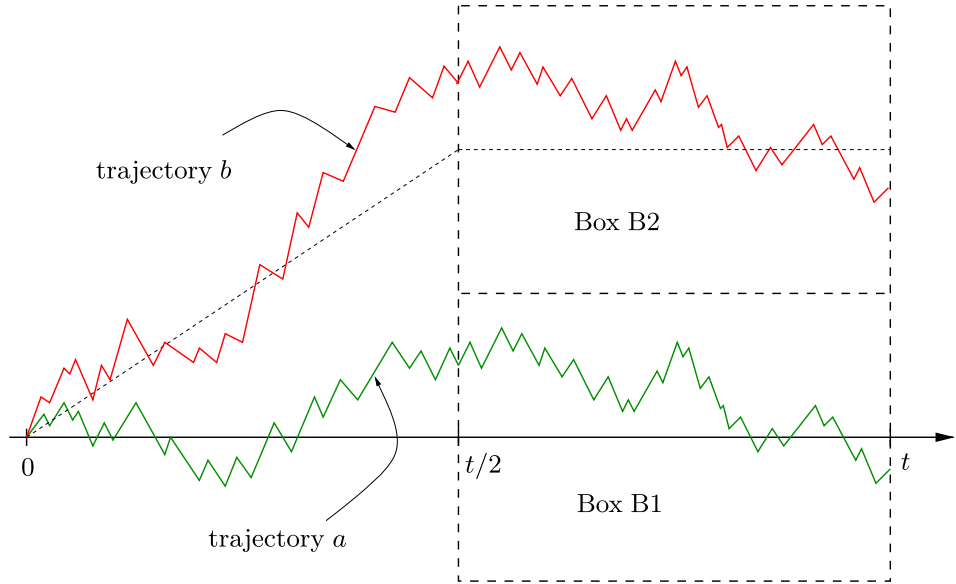}

\caption{A representation of the two events for which we want to
compare weights, and two typical trajectories for each event ($a$ and
$b$).} \label{superdiff}
\end{figure}

We look at the weight under $\mu_t$ of the trajectories $(B_s)_{s\in
[0,t]}$ that stay within a box of width $t^{\alpha}$ centered on the
origin for $s\in[t/2,t]$ (box $B1$) and compare this with the
weight of
the trajectories that spend all the time $s\in[t/2,t]$ in another box
of the same width (see Figure \ref{superdiff}, trajectories $a$ and
$b$). The entropic cost for one trajectory to reach the box $B2$ and
stay there is equal to $\log P\{B \mbox{ stays in } B2\}
\sim t^{2\alpha-1}$.

%

In order to estimate the energy variation between the two boxes, we
look at $\gO_i,$ the sum of all the increments of $\omega$ in the box $Bi$
($i=1,2$):
\[
\gO_i :=\int_{Bi} \omega(\dd s, x)\,\dd x.
\]

The $\gO_i$ are Gaussian variables that are identically distributed,
with variance~$\approx t^{\alpha(2-\theta)+1}$. Therefore, in each box,
the empirical mean of the increment of $\omega$ by a unit of time in $Bi$
($\gO_i/|Bi|$) is Gaussian with variance $\approx t^{ -\alpha\theta-1}$
[the typical fluctuations are of order $t^{ (-\alpha\theta-1)/2}$].
Multiplying this by the length of the box ($t/2$), we get that the
empirical  mean of the energy for paths in a given box has typical fluctuations of
order $t^{(1-\alpha\theta)/2}$.

Therefore, if $\alpha<\frac{3}{4-\theta}$, with positive probability,
the energy gain for going into the box $B2$ is at least
$t^{(1-\alpha)/2}$ and is bigger than the entropy cost
$t^{2\alpha-1}$, so the trajectory is less likely to stay in the box
$B1$ than in the box $B2$.

To make this argument rigorous, we have to:
\begin{longlist}[$\bullet$]
\item[$\bullet$] use Girsanov path transforms to make the argument about the
entropy work;
\item[$\bullet$] use a measure coupling argument to make the energy comparison rigorous;
\item[$\bullet$] make the comparison with more than two boxes.
\end{longlist}

In this section, for practical reasons, we work with $\| \cdot
\|_{\infty}$, the $l_\infty$-norm in $\R^d$, and not with the Euclidean
norm.

\subsection[Proof of Theorem 1.6]{Proof of Theorem \protect\ref{supdiffff}}

Let $N$ be some even integer, $\alpha:=\frac{3}{4+\theta}$. For
$k\in
\N,$ define
\[
\gL_k:=[t/2,t]\times\biggl[\frac{(2k-1)t^{\alpha}}{N^2},\frac
{(2k+1)t^{\alpha}}{N^2}\biggr]
\times\biggl[-\frac{t^{\alpha}}{N^2},\frac{t^{\alpha}}{N^2}\biggr]^{d-1}.
\]
Define
\[
Z_t^{(k)}:= P\biggl[\exp\biggl(\gb\int_0^t \omega(\dd
s,B_s)\biggr)\ind_{\{(s,B_s)\in\gL_k, \forall s\in[t/2,t]\}}\biggr].
\]
The proof can be decomposed in two steps. The first step (the next
lemma) has been strongly inspired by the work of Petermann \cite
{cfP}, Lemma 2, and gives a rigorous method to bound from above the entropic
cost for reaching a region $t^{\alpha}$ away from the origin.

\begin{lemma}\label{th:Q}
With probability greater than $1-1/N$, we have
\[
\sum_{k\in\{-N,\dots,N\}\setminus\{0\}} Z_t^{(k)}\ge
\exp\biggl(-\frac{8}{N^2}t^{2\alpha-1}\biggr)Z_t^{(0)}.
\]
\end{lemma}
\begin{pf}
We use the transformation $h_k\dvtx\bbR_+\times\bbR^d\rightarrow
\bbR^d
$ which transforms a path contributing to $Z_t^{(0)}$ into a path
contributing to $Z_t^{(k)}$.
\[
h_k\dvtx(s,x)\mapsto x+\bigl((2s/t)\wedge1\bigr)\frac
{2k}{N^2}t^{\alpha}\mathbf{e}_1,
\]
where $\mathbf{e}_1$ is the vector $(1,0,\dots,0)$ in $\bbR^d,$ and we
define
\[
\barr Z_t^{(k)}:= P\exp\biggl(\gb\int_0^t \omega(\dd s,
h_k(s,B_s))\biggr)\ind_{\{(s,B_s)\in A_0, \forall s\in[t/2,t]\}}.
\]
One can check that $( \barr Z_t^{(k)}
)_{k\in\{-N/2,\dots,N/2\}}$ is a family of identically
distributed random variables. Moreover, elementary reasoning gives us
that there exists an integer $k_0,$ with $|k_0|\le N/2,$ such that
\[
Q\Bigl\{\barr Z_t^{(k_0)}= \max_{k\in\{-N/2,\dots,N/2\}} \barr
Z_t^{(k)}\Bigr\}\le\frac{1}{N+1}.
\]
From this, one infers, by translation invariance of the environment,
that
\[
Q\Bigl\{\barr Z_t^{(0)}\ge \max_{k\in\{-N/2-k_0,\dots,N/2-k_0\}
}\barr Z_t^{(k)}\Bigr\}\le\frac{1}{N+1}.
\]
The final step to get the result is to compare $\barr Z^{(k)}_t$ with
$Z_t^{(k)}$. This can be done by using a Girsanov path transform:
\begin{eqnarray*}
&&Z_t^{(k)}=P \exp\biggl(-\frac{4k}{N^2}t^{\alpha
-1}B^{(1)}_{t/2}-\frac{4k^2}{N^4}t^{2\alpha-1}+\gb\int_0^t \omega
(\dd s ,h_k(s,B_s))\biggr)\\
&&\quad\hphantom{Z_t^{(k)}}{}\times\ind_{\{(s,B_s)\in A_0, \forall
s\in[t/2,t]\}}\\
&&\hphantom{Z_t^{(k)}}\ge\exp\biggl(-\frac
{4(k^2+|k|)}{N^4}t^{2\alpha-1}\biggr)\barr Z_t^{(k)}\ge\exp\biggl
(-\frac{8}{N^2}t^{2\alpha-1}\biggr)\barr Z_t^{(k)},
\end{eqnarray*}
where $B^{(1)}$ is the first coordinate of the Brownian motion.
\end{pf}

\begin{rem}
If one chooses $\alpha=1/2,$ this lemma shows that in any dimension,
for any correlation function $Q$, the directed polymer is at least
diffusive, that is, that the typical distance to the origin is of order
at least $\sqrt{N}$.
\end{rem}

For the rest of the proof, the idea which is used differs substantially
from the one used in \cite{cfP} and then adapted in \cite{cfBTV}.
Instead of using purely Gaussian tools and working with covariance
matrices, we use changes of measure that are similar to those used in
the previous section. This shortens the proof considerably and makes it
less technical and more intuitive. Moreover, it highlights the fact
that the proof could be adapted to a non-Gaussian context, for
example, the model of Brownian polymer in a Poissonian environment
studied by Comets and Yoshida \cite{cfCYcmp}.

We set $T:=t^{\alpha}N^{-3}$ and define
\[
\Omega:=\frac{\int_{[-T,T]^d}\int_{t/2}^{t}\omega(\dd s,x)\,\dd
x}{\sqrt{t/2 \int_{[-T,T]^d\times[-T,T]^d}Q(x-y)\,\dd x\,\dd y}},
\]
which is a standard centered Gaussian variable. We define the
probability measure $\bP_0$ on the environment by its Radon--Nikodym
derivative with respect to $\bP$:
\[
\frac{\dd\bP_0}{\dd\bP}(\omega):=\exp(-\gO-1/2).
\]

The probability $\bP_0$ has two very important characteristics. It is
not very different from $\bP$ (see the next lemma) and it makes the
environment less favorable in the box $[-t/2,t]\times[-T,T]^d$ so that,
under $\bP_0,$ the trajectory will be less likely to stay in that box.

\begin{lemma}\label{petitlem}
Letting $A$ be any event, we have
\[
\bP(A)\le\sqrt{e\bP_0(A)}.
\]
\end{lemma}
\begin{pf}
This is a simple application of the H\"older inequality:
\[
\bP(A)=\bP_0 \biggl[\frac{\dd\bP}{\dd\bP_0} \ind_{A}\biggr
]\le\sqrt{\bP\biggl[\frac{\dd\bP}{\dd\bP_0}\biggr]}\sqrt{\bP
_0 (A)}.
\]
\upqed
\end{pf}

Now, our aim is to show that under $\bP_0$, the probability that the
walk stays in $[-t^{\alpha}/N^3,t^{\alpha}/N^3]^d$ is small and then to
use the above lemma to conclude. We use Lemma \ref{gaussest} to define,
on the same space, two environments with laws $\bP$ and $\bP_0$.
Indeed, if $\omega$ has distribution $\bP,$ then $\hatt\omega$
defined by
%
\begin{eqnarray}\label{coopling}
\hatt\omega(0,x)&:=& 0 \qquad\forall x\in\R^d,\nonumber\\
\hatt\omega( \dd s,x)&:=& \omega(\dd s,x)-\bP[\gO\omega(\dd
s,x)]\\
&=& \omega(\dd s,x)-\frac{\ind_{s\in[t/2,t]}\int_{[-T,T]^d}Q(x-y)\,
\dd y}{\sqrt{t/2 \int_{[-T,T]^d\times[-T,T]^d}Q(x-y)\,\dd x\,\dd
y}}\,\dd s\nonumber
\end{eqnarray}
has distribution $\bP_0$. We define
\begin{eqnarray*}
X_t&:=&P\biggl[ \exp\biggl(\gb\int_0^t \omega(\dd s,B_s)\biggr
)\ind_{\{B_s\in[-T,T]^d,\forall s\in[t/2,t]\}}\biggr],\\
\hatt X_t&:=& P \biggl[\exp\biggl(\gb\int_{0}^{t} \hatt\omega
(\dd s,B_s)\biggr)\ind_{\{B_s\in[-T,T]^d,\forall s\in[t/2,t]\}
}\biggr],\\
\hatt Z_t^{(k)}&:=& P\biggl[\exp\biggl(\gb\int_0^t \hatt\omega
(\dd s,B_s)\biggr)\ind_{\{(s,B_s)\in\gL_k,\forall s\in[t/2,t]\}
}\biggr].
\end{eqnarray*}
From this definition
%
\begin{equation}\label{truc}
\mu^{\gb,\omega}_t\Bigl\{\max_{s\in[0,t]}\|B_s\|_{\infty} \le
t^{\alpha}/N^3\Bigr\}\le\frac{X_t}{\sum_{k\in\{-N,\dots,N\}
\setminus\{0\}} Z_t^{(k)}}.
\end{equation}
And, for any $x$,
%
\begin{eqnarray}\label{der22}
&&\bP_0\Bigl[ \mu_t^{\gb,\omega}\Bigl\{\max_{s\in[0,t]} \|B_s\|
_{\infty}\le t^{\alpha}/N^3\Bigr\}\le x \Bigr]\nonumber\\
&&\qquad= \bP\Bigl[ \mu_t^{\gb,\hatt\omega}\Bigl\{\max_{s\in
[0,t]}\|B_s\|_{\infty}\le t^{\alpha}/N^3\Bigr\} \le x \Bigr]\\
&&\qquad\le\bP\biggl[\frac{\hatt X_t}{\sum_{k\in\{-N,\ldots,N\}
\setminus\{0\}} \hatt Z_t^{(k)}}\le
x\biggr].\nonumber
\end{eqnarray}
We will now use measure coupling to bound the right-hand side of the
above.

\begin{lemma}
For $N$ large enough and $t \gg N^{3/\alpha}$, we have
\[
\frac{\hatt X_t}{\sum_{k\in\{-N,\dots,N\}\setminus\{0\}} \hatt
Z_t^{(k)}}\le\frac{X_t}{\sum_{k\in\{-N,\dots,N\}\setminus\{0\}}
Z_t^{(k)}}\exp(-C_4 t^{1/2}T^{-\theta/2}).
\]
\end{lemma}
\begin{pf}
From the definition of $X_t$ and $\hatt X_t$, we have
\[
\log(\hatt X_t/X_t)\le \gb\sup_{\{B:\|B_s\|_{\infty} \le T,\forall
s \in[t/2,t]\}} \int_{0}^t \bigl(\hatt\omega(\dd s,B_s)-\omega
(\dd s, B_s)\bigr),
\]
where the ``sup'' is to be understood as the essential supremum under
the Wiener measure $P$. It follows from the coupling construction
\eqref{coopling} that for the trajectories staying within $[-T,T]^d$ in
the time interval $[t/2,t],$ we have
\begin{eqnarray*}
\int_{0}^t \bigl(\hatt\omega(\dd s,B_s)-\omega(\dd s, B_s)\bigr)&=&-\bP\biggl[\int_{t/2}^t \omega(\dd s,B_s)\gO\biggr]\\
&=&-\frac{\ind_{[-T,T]^d}\int_{t/2}^t Q(B_s-y)\,\dd s\,\dd y
}{\sqrt{t/2 \int_{[-T,T]^d\times[-T,T]^d}Q(x-y)\,\dd x\,\dd y}}\,\dd t\\
&\le&-\frac{\sqrt{t/2}\inf_{x\in[-T,T]^d} \int_{[-T,T]^d}Q(x-y)\,\dd y}{\sqrt{\int_{[-T,T]^d\times[-T,T]^d}Q(x-y)\,\dd x\,\dd y}}
\end{eqnarray*}
so that
%
\begin{equation}\label{xxxest}
\log(\hatt X_t/X_t)\le - \frac{\gb\sqrt{t/2}\inf_{x\in[-T,T]^d}
\int_{[-T,T]^d} Q(x-y)\,\dd y}{\int_{[-T,T]^d\times[-T,T]^d}Q(x-y)\,
\dd x\,\dd y}.
\end{equation}

Performing a similar computation, one gets that
%
\begin{eqnarray}\label{zzzest}
&&\log\frac{\sum_{k\in\{-N,\dots,N\}\setminus\{0\}}
Z_t^{(k)}}{\sum_{k\in\{-N,\dots,N\}\setminus\{0\}} \hatt
Z_t^{(k)}}\nonumber\\ [-8pt]\\ [-8pt]
&&\qquad\le \frac{\gb\sqrt{t/2}\max_{\|x\|_{\infty}\ge t^{\alpha
}/N^2} \int_{[-T,T]^d} Q(x-y)\,\dd y}{\sqrt{\int_{[-T,T]^d\times
[-T,T]^d}Q(x-y)\,\dd x\,\dd
y}}.\nonumber
\end{eqnarray}
It remains to give an estimate for the right-hand sides of
\eqref{xxxest} and \eqref{zzzest}.

First, we remark that
%
\begin{equation}\label{joli1}
\int_{[-T,T]^d\times[-T,T]^d} Q(x-y)\,\dd x\,\dd y\asymp T^{2d-\theta}.
\end{equation}
For \eqref{xxxest}, we note that one can find a constant $C_5$ such
that for all $x,y \in[-T,T]^d$, $Q(x-y)\ge C_5 T^{-\theta}$, so that
%
\begin{equation}\label{joli2}
\inf_{x\in[-T,T]^d} \int_{[-T,T]^d} Q(x-y)\,\dd y \ge C_5
T^{d-\theta}.
\end{equation}
For \eqref{zzzest}, we note that one can find a constant $C_6$ such
that for all $y \in[-T,T]^d$, $\|x\|_{\infty}\ge t^{\alpha}/N^{2}$,
$Q(x-y)\le C_6 (t^{\alpha}/N^2)^{-\theta}=C_6 (N T)^{-\theta}$, so that
%
\begin{equation}\label{joli3}
\max_{\|x\|_{\infty}\ge t^{\alpha}/N^2} \int_{[-T,T]^d} Q(x-y)\,\dd
y \le C_6 N^{-\theta} T^{d-\theta}.
\end{equation}
For $N$ large enough, the term above will be dominated by $C_5
T^{d-\theta}$ so that, combining \eqref{joli1}, \eqref{joli2} and
\eqref{joli3}, one gets that there exists $C_4$ such that
\[
\log\frac{\hatt X_t\sum_{k\in\{-N,\dots,N\}\setminus\{0\}}
Z_t^{(k)}}{X_t \sum_{k\in\{-N,\dots,N\}\setminus\{0\}} \hatt
Z_t^{(k)}}\le-C_4 t^{1/2} T^{-\theta/2}.
\]
\upqed
\end{pf}

Now, the preceding result, together with Lemma \ref{th:Q}, ensures
that, with probability larger than $1-1/N,$ we have
%
\begin{eqnarray}\label{der11}
\qquad \frac{\hatt X_t}{\sum_{k\in\{-N,\dots,N\}\setminus\{0\}} \hatt
Z_t^{(k)}}&\le& \frac{Z_t^{(0)}}{\sum_{k\in\{-N,\dots,N\}\setminus
\{0\}}
Z_t^{(k)}}\exp(-t^{1/2}T^{-\theta/2})\nonumber\\ [-8pt]\\ [-8pt]
&\le&\exp\biggl(\frac{8}{N^2}t^{2\alpha-1}-C_4t^{1/2}T^{-\theta
/2}\biggr).\nonumber
\end{eqnarray}
We can bound the term in the exponential on the right-hand side when
$t$ is large enough:
%
\begin{eqnarray}\label{der12}
\frac{8}{N^2}t^{2\alpha-1}-C_4t^{1/2}T^{-\theta/2}&=&\biggl[\frac{8}{N^2}-C_4
N^{3\theta/2}\biggr]t^{(2-\theta)/(4+\theta)}\nonumber\\ [-8pt]\\ [-8pt]
&\le&-t^{(2-\theta)/(4+\theta)}.\nonumber
\end{eqnarray}
We now combine all of the elements of our reasoning. Equations
\eqref{der11} and \eqref{der12} combined with \eqref{der22} give us
\[
\bP_0\Bigl\{\mu^{\gb,\omega}_t\Bigl\{\max_{0\le s\le t} \|B_s\|
_{\infty} \le t^{\alpha}/N^3\Bigr\}\ge\exp\bigl(-t^{(2-\theta
)/(4+\theta)}\bigr)\Bigr\}\le\frac{1}{N}.
\]
Lemma \ref{petitlem} allows us to get from this
\[
\bP\Bigl\{\mu_t^{\gb,\omega}\Bigl\{\max_{0\le s\le t} \|B_s\|
_{\infty} \le t^{\alpha}/N^3\Bigr\}\ge\exp\bigl(-t^{(2-\theta
)/(4+\theta)}\bigr)\Bigr\}\le\sqrt{e/N}
\]
so that
\[
\bP\mu_t^{\gb,\omega} \{\max\|B_s\|_{\infty} \le t^{\alpha}/N^3\}
\le\exp\bigl(-t^{(2-\theta)/(4+\theta)}\bigr)+\sqrt{e/N}.
\]
We obtain the desired result by choosing $N$ large enough.\qed

\begin{rem}
The above proof is valid for $\theta<2$, $d\ge2$ and also for $\theta
<1$, $d=1$.
For the case $d=1$, $\theta>1$ [or, more generally, for $d=1$, $Q\ge
0$, $\int_{\R}Q(x)\,\dd x<\infty$], one has to choose $\alpha=3/5$ and
replace \eqref{joli1} by
\[
\int_{[-T,T]\times[-T,T]} Q(x-y)\,\dd x\,\dd y\asymp T
\]
and \eqref{joli2}, \eqref{joli3} by
\begin{eqnarray*}
\inf_{x\in[-T,T]} \int_{[-T,T]} Q(x-y)\,\dd y&\ge& C_5,\\
\max_{x\ge T^{\alpha}/N^2} \int_{[-T,T]} Q(x-y)\,\dd y&\le& C_5/2.
\end{eqnarray*}
\end{rem}

\subsection[Proof of Theorem 1.11]{Proof of Theorem \protect\ref{abcdef}}

Let $\alpha$ be as in the assumption of the theorem. We define
$T:=t^{\alpha}$ and
\[
\check Z_t:= \bP\bigl[\exp(\gb H_t(B))\ind_{\{\max_{s\in[0,t]} \|
B_s\|_{\infty}\le t^{\alpha}\}}\bigr].
\]
We will show that the typical fluctuations of $\log\check Z_t$ are
large and then that those of $\log Z_t$ are also large because $Z_t$
and $\check Z_t$ are close to each other with large probability (in
this respect, we need to do more than simply bound the variance of
$\log\check Z_t$). Define $\gO$ by

\[
\gO:=\frac{\int_{[0,t]\times[-T,T]^d}\omega(\dd s,x)\,\dd x}{\sqrt
{t\int_{[-T,T]^{2d}}Q(x-y)\,\dd x\,\dd y}}.
\]
Under $\bP$, $\gO$ is a standard Gaussian. We define the probability
$\bP_0$ by
\[
\dd\bP_0:=\exp(\gO-1/2)\,\dd\bP.
\]
It follows from its definition that the distribution of $\log\check
Z_t$ under $\bP$ has no atom, so one can define $x_0$ (depending on
$\gb$ and $t$) such that
\[
\bP[ \log\check Z_t\le x_0]=e^{-2}.
\]

We use Lemma \ref{gaussest} to perform a measure coupling as before: if
$\omega$ has distribution $\bP$, then $\hatt\omega$ defined by
\begin{eqnarray*}
\hatt\omega(0,x)&:=&0 \qquad\forall x\in\R^d,\\
\hatt\omega( \dd s,x)&:=& \omega(\dd s,x)+\bP[\gO\omega(\dd
s,x)]\\
&=& \omega(\dd s,x)+\frac{\int_{[-T,T]^d}Q(x-y)\,\dd y}{\sqrt{t\int_{[-T,T]^{2d}}Q(x-y)\,\dd x\,\dd y}}\,\dd s
\end{eqnarray*}
has distribution $\bP_0$.

One defines
\[
\hatt Z_t:= P\biggl[\exp\biggl(\gb\int_0^t \hatt\omega(\dd s,B_s)\biggr)\ind_{\{ \max_{s\in[0,t]} \| B_s \|_{\infty} \le T
\}}\biggr].
\]

For all paths $B$ such that $\max_{s\in[0,t]} \| B_s \|_{\infty} \le
T$, we have
\begin{eqnarray*}
\int_0^t \hatt\omega(\dd s, B_s)-\int_0^t \omega(\dd s,
B_s)&=&\int_0^t \frac{\int_{[-T,T]^d}Q(B_s-y)\,\dd y}{\sqrt{t \int
_{[-T,T]^2d}Q(x-y)\,\dd x\,\dd y}}\,\dd s\\
&\ge& \sqrt{t} \frac{\min_{\|x\|_{\infty}\le T}\int
_{[-T,T]^d}Q(x-y)\,\dd
y}{\sqrt{\int_{[-T,T]^{2d}}Q(x-y)\,\dd x\,\dd y}}.
\end{eqnarray*}
One can check that there exists a constant $c$ depending only on $Q$
such that for all~$t$,
\[
\frac{\min_{\|x\|_{\infty}\le T}\int_{[-T,T]^d}Q(x-y)\,\dd y}{\sqrt
{\int_{[-T,T]^{2d}}Q(x-y)\,\dd x\,\dd y}} \ge c T^{-(\theta\wedge d)/2}.
\]
Therefore, one has $\bP$ almost surely
\[
\label{trucmachin}
\log\hatt Z_t\ge\log\check Z_t+c\gb t^{(1-\alpha(\theta\wedge d))/2}.
\]
We use Lemma \ref{petitlem} (which is still valid in our case, even if
the change of measure is different) to see that
\[
\bP\bigl\{\log\check Z_t\le x_0+ c\gb t^{(1-\alpha(\theta
\wedge d))/2}\bigr\}
\le\sqrt{e \bP_0\bigl\{ \log\check Z_t\le x_0+ c\gb
t^{(1-\alpha(\theta\wedge d))/2}\bigr\}}.
\]
Moreover,
\begin{eqnarray*}
\bP_0\bigl\{ \log\check Z_t\le x_0+ c\gb t^{(1-\alpha(\theta
\wedge d))/2}\bigr\}&=&\bP\bigl\{\log\hatt Z_t\le x_0+ c\gb t^{(1-\alpha(\theta
\wedge d))/2}\bigr\}\\
&\le&\bP\{\log\check Z_t \le x_0\}=e^{-2}.
\end{eqnarray*}
The inequality is given by \eqref{trucmachin}. So, combining
everything, we have
%
\begin{eqnarray}\label{fghuio}
\bP\bigl\{\log\check Z_t\ge x_0+ c\gb t^{(1-\alpha(\theta\wedge
d))/2}\bigr\}&\ge&
1-e^{-1/2},\nonumber\\ [-8pt]\\ [-8pt]
\bP\{\log\check Z_t\le x_0\}&=& e^{-2}.\nonumber
\end{eqnarray}
It is enough to prove that the variance of $\log\check Z_t$ diverges
with the correct rate. Slightly more work is required to prove the same
for $\log\check Z_t$. Recall that
\[
\mu^{\gb,\omega}_t\Bigl\{ \max_{s\in[0,t]} \| B_s \|_{\infty}
\le T\Bigr\}= \check Z_t/Z_t.
\]
The assumption on $\alpha$ gives that, for $t$ large enough,
\[
\bP\{ (\check Z_t/Z_t)\le1/2\}\le1/100.
\]
By a union bound, we have
\[
\bP\{\log\check Z_t\le x_0 \}\le\bP\{\log Z_t\le x_0+\log2 \}+
\bP\{ \check Z_t \le Z_t/2 \}.
\]

Combining this with \eqref{fghuio} (and the trivial bound $Z_t\ge
\check Z_t$) gives that for $t$ large enough,
\begin{eqnarray*}
\bP\bigl\{\log Z_t\ge x_0+ c\gb t^{(1-\alpha(\theta\wedge
d))/2}\bigr\}&\ge& 1-e^{-1/2},\\
\bP\bigl\{\log Z_t\le x_0+\log2 \bigr\}&\ge& e^{-2}-1/100.
\end{eqnarray*}
This implies that
\[
\hspace*{30pt}\var\log Z_t\ge (e^{-2}-1/100)(1-e^{-1/2})\bigl(c\gb t^{(1-\alpha
(\theta\wedge d))/2}-\log2\bigr)^2.\hspace*{30pt}\qed
\]

\subsection[Proof of Proposition 1.9]{Proof of Proposition \protect\ref{mmeejj}}

To prove this result, we will follow the method of M\'ejane
\cite{cfM}. First, we need to use a concentration result. We prove it
using stochastic calculus.
\begin{lemma}\label{concent}
Let $f$ be a nonnegative function on $\R^d$ and $r\in[0,t]$. We
define $\tilde Z_t:=P[f(B_r)\exp(\gb\int_0^t \omega(\dd s, B_s))]$.
Then, for all $x> 0$,
\[
\bP\bigl\{ |\log\tilde Z_t -\bP[ \log\tilde Z_t ]|\ge
x\sqrt{t} \bigr\}\le2\exp\biggl(-\frac{x^2}{2\gb^2}\biggr).
\]
\end{lemma}

\begin{pf}
Let $\mathcal F=(\mathcal F_t)_{t\ge0}:=(\sigma\{\omega(s,x),
s\in[0, t], x\in\R^d\})_{t\ge0}$ be the natural
filtration associated with the environment. We consider the following
continuous martingale with respect to the filtration $\mathcal F$:
%
\begin{equation}\label{mudef}
(M_u:=\bP[\log\tilde Z_t \big| \mathcal F_u]-\bP[\log\tilde Z_t
])_{u\in[0,t]}.
\end{equation}
We have $M_0=0$ and the result to prove becomes
\[
\bP\bigl\{|M_t|\ge x\sqrt{t}\bigr\}\le2 \exp\biggl(-\frac
{x^2}{2\gb^2}\biggr).
\]

The proof uses a classical result on concentration for martingales.
\begin{lemma}\label{martin}
If $(M_u)_{u\ge0}$ is a continuous martingale (with associated law
$\bP$) starting from $0$ with finite exponential moments of all
orders, then, for all $(x,y)\in\R_+^2$ and $u\ge0,$ we have
\[
\bP\{M_{u}\ge x\ll M\rr_{u}+y\}\le\exp(-2xy).
\]
\end{lemma}
\begin{pf}
We have
\begin{eqnarray*}
\bP\{M_{u}\ge x\ll M\rr_{u}+y\}&=&\bP\{\exp(2xM_{u}-2x^2\ll M\rr
_u-2xy)\ge1\}\\
&\le& \bP[\exp(2xM_{u}-2x^2\ll M\rr_u-2xy)]=\exp(-2xy),
\end{eqnarray*}
where we have just used the fact that for any given $x$,
$\exp(xM_{u}-2x^2\ll M\rr_u)$ is a martingale.
\end{pf}

To use the previous lemma, we have to compute the bracket of the
martingale $M$ defined in $\eqref{mudef}$. One can compute it
explicitly, but the form of the result is rather complicated so that we
have to introduce several items of notation before giving a formula.
One defines the probability measure $(\tilde\mu_t)$ (depending on
$\omega$, $\gb$, $t$ and $r$) by giving its Radon--Nikodym derivative with
respect to the Wiener measure:

\[
\frac{\dd\tilde\mu_t(\dd x)}{\dd P}(B):=\frac{1}{\tilde
Z_t}f(B_r)\exp\biggl(\gb\int_0^t \omega(\dd s, B_s)\biggr)
\]
[it is the polymer measure where the paths have been reweighted by
$f(B_r)$]. One defines $\nu_{t,u}$ to be the (random) probability
measure on $\R$ defined by

\[
\nu_{t,u}(\dd x):=\bP[\tilde\mu_t(B_u\in\dd x)\vert\mathcal F_u].
\]
For the martingale $M$ defined by $\eqref{mudef}$, we have
%
\begin{equation}\label{baraket}
\ll M\rr_t= \gb^2 \int_0^t \int_{\bbR^2} \nu^{\otimes
2}_{t,u}(\dd x\,\dd y)Q(x-y).
\end{equation}
An easy consequence of this is that
\[
\ll M\rr_t\le\gb^2 t \qquad\mbox{almost surely}.
\]

From this, we infer that
\[
\bP\bigl\{M_t\ge x\sqrt{t}\bigr\}\le\bP\biggl\{M_t\ge\frac
{x}{2\gb^2\sqrt{t}}\ll M\rr_t+\frac{\sqrt{t}x}{2}\biggr\}\le\exp
(-x^2\gb^{-2}/2),
\]
where the last inequality is obtained by applying Lemma \ref{martin}.
Carrying out the same computation for the martingale $-M$ gives the
desired result.
\end{pf}

We now turn to the proof of the result. Let $\alpha>3/4$ be fixed. Let
$B^{(1)}_t$ be the first coordinate of $B_t$ in $\R^d$. By the Markov
inequality, we have, for every $\lambda>0,$
%
\begin{equation}\label{mutb}
\mu_t\bigl\{ B^{(1)}_r\ge a\bigr\}\le e^{-\lambda a+r\lambda
^2/2}\mu_t\bigl(\exp\bigl(\lambda B^{(1)}_r-r\lambda^2/2\bigr
)\bigr).
\end{equation}
We use Girsanov's formula:
%
\begin{eqnarray}\label{girsasa}
\mu_t\bigl(\exp\bigl(-\lambda B^{(1)}_r-r\lambda^2/2\bigr)\bigr
)&=&\frac{P[\exp(\lambda B^{(1)}_r-r\lambda^2/2+\gb\int
_{0}^{t}\omega(\dd s, B_s))]}{P[\exp(\gb\int_{0}^t\omega(\dd s,
B_s))]}\nonumber\\ [-8pt]\\ [-8pt]
&=&\frac{P[\exp(\gb\int_{0}^t\omega(\dd s,h^{\lambda
,r}(s,B_s)))]}{P[\exp(\gb\int_0^t\omega(\dd s,B_s))]},\nonumber
\end{eqnarray}
where $h^{\lambda,r}$ is the function from $\R_+\times\R^d$ to $\R^d$
defined by
\[
h^{\lambda,r}(s,x):=x+\lambda(r\wedge s){\bf e}_1
\]
and ${\bf e}_1$ is the vector $(1,0,\dots,0)$ in $\R^d$. By translation
invariance, the environment $(\omega(s,
h^{\lambda,r}(s,x)))_{s\in[0,t],x\in\R^d}$ has the \vspace{1pt}same law as
$(\omega(s,x))_{s\in[0,t],x\in\R^d}$ and, therefore, we get,
from the last line of \eqref{girsasa}, that
\[
\bP\bigl[\log\mu_t\bigl(\exp\bigl(-\lambda B^{(1)}_r-r\lambda
^2/2\bigr)\bigr)\bigr]=0.
\]
Substituting this into \eqref{mutb}, we get
\[
\bP\log\mu_t\bigl\{ B^{(1)}_r\ge a\bigr\}\le-\lambda a+r\lambda^2/2.
\]
As $\lambda$ is arbitrary, we can take the minimum over $\lambda$ for the
right-hand side to get $-\frac{a^2}{2r}$. We use the result for
$a=t^{\alpha}$ to get
\[
\bP\log\mu_t\bigl\{ B^{(1)}_r\ge t^{\alpha}\bigr\}\le-t^{2\alpha}/2r.
\]
Using Lemma \ref{concent} with $f(y):=\ind_{\{y\ge a\}}$ and
$f\equiv1$, for $x=t^{\gep}$ with $\gep<4\alpha-3,$ one
gets
\[
\bP\biggl\{\log\mu_t\bigl\{ B^{(1)}_r\ge t^{\alpha} \bigr\}\le
-\frac{t^{2\alpha}}{2r}+2
t^{(1+2\gep)/2}\biggr\}\le4\exp\biggl(-\frac{t^{2\gep}}{2\gb
^2}\biggr).
\]
For $t$ sufficiently large, we have, for all $r\le t$,
\[
-\frac{t^{2\alpha}}{2r}+2 t^{(1+2\gep)/2}\le
-t^{1/2}.
\]

We can get this inequality for $B^{(i)}$ and $-B^{(i)}$ for any
$i\in\{1,\dots,d\}$. Combining all of these results, we get
\[
\bP\bigl\{\mu_t\{ \|B_r\|_{\infty}\ge t^{\alpha}\}\le2d\exp
(-t^{1/2})\bigr\}\le8d\exp\biggl(-\frac{t^{2\gep}}{2\gb^2}\biggr).
\]
Using the above inequality for all $r\in\{1,2,\dots,\lfloor t
\rfloor\}$, we obtain that
%
\begin{eqnarray}\label{partie1}
&&\bP\Bigl\{\mu_t\Bigl\{ \max_{ r\in\{1,2,\dots,\lfloor t
\rfloor\}} \|B_r\|_{\infty}\ge t^{\alpha} \Bigr\}\le2dt
\exp(-t^{1/2})\Bigr\}\nonumber\\ [-8pt]\\ [-8pt]
&&\qquad\le8dt\exp\biggl(-\frac{t^{2\gep}}{2\gb^2}\biggr
).\nonumber
\end{eqnarray}
To complete the proof, we need to control the term
\[
\bP\mu_t\Bigl\{\mathop{\max_{s\in[n,n+1]}}_{n\in\{0,\dots
,\lfloor
t\rfloor\}} \|B_s-B_n\|_{\infty}\ge t^{\alpha} \Bigr\}.
\]
One computes that
%
\begin{eqnarray}\label{trucpetit}
&&\bP\Bigl[W_t \mu_t\Bigl\{\mathop{\max_{s\in[n,n+1]}}_{n\in\{
0,\dots,\lfloor t\rfloor\}}\|B_s-B_n\|_{\infty}\ge
t^{\alpha}\Bigr\}\Bigr]\nonumber\\ [-8pt]\\ [-8pt]
&&\qquad=P\Bigl\{\mathop{\max_{s\in[n,n+1]}}_{n\in\{0,\dots
,\lfloor t\rfloor\}} \|B_s-B_n\|_{\infty}\ge t^{\alpha}\Bigr\}
\le4d(t+1)\exp(-t^{2\alpha}/2).\nonumber
\end{eqnarray}
Lemma \ref{concent} applied to $f\equiv1$ and $x=\gb\sqrt{t}$ gives us
that
%
\begin{equation}\label{trucpetit2}
\bP\bigl\{W_t\le\exp\bigl(t\bigl(p(\gb)-\gb\bigr)\bigr)\bigr
\}\le\exp(-t/2).
\end{equation}
If $A$ is an event and $x>0$ are given, then
\begin{eqnarray*}
\bP[\mu_t(A)]&=&\bP\bigl[\mu_t(A)\ind_{\{W_t\le x\}}\bigr]+\bP
\bigl[\mu_t(A)\ind_{\{W_t> x\}}\bigr]\\
&\le& \bP\{W_t\le x\}+x^{-1}\bP[W_t\mu_t(A)].
\end{eqnarray*}
Therefore, we have, from \eqref{trucpetit2} and \eqref{trucpetit}, that
%
\begin{eqnarray}\label{partie2}
\qquad\bP\mu_t\Bigl\{\mathop{\max_{s\in[n,n+1]}}_{n\in\{0,\dots,\lfloor t\rfloor\}} \|B_s-B_n\|_{\infty}\ge t^{\alpha} \Bigr\}&\le&
e^{-t/2}+4d(t+1)e^{t(\gb-p(\gb))-t^{2\alpha}/2}\nonumber\\[-8pt]\\ [-8pt]
&\le& 2 e^{-t/2}.\nonumber
\end{eqnarray}
Combining \eqref{partie1} and \eqref{partie2}, we get (for $t$ large
enough)
\begin{eqnarray*}
\hspace*{43pt}\bP\mu_t\Bigl\{\max_{s\in[0,t]} \|B_s\|_{\infty}\le2t^{\alpha}
\Bigl\}&\le& 2 \exp(-t/2)+8dt\exp(-t^{2\gep}/\gb^2)\\
&\le& \exp(-t^{\gep}).\hspace*{134.5pt}\qed
\end{eqnarray*}

\section*{Acknowledgments} The author is very grateful to Markus Petermann
for sending him a copy of his unpublished work, and to Giambattista
Giacomin for much enlightening discussion and advice.

%

\printaddresses

\end{document}